\documentclass[EJP]{ejpecp} 

\usepackage[utf8]{inputenc}

\usepackage{tikz}

\SHORTTITLE{An extension of the Ising-CW model of SOC with a threshold on the interaction range}

\TITLE{An extension of the Ising-Curie-Weiss model of self-organized criticality with a threshold on the interaction range}

\AUTHORS{
  Nicolas~Forien
  \footnote{
  Aix Marseille Univ, CNRS, Centrale Marseille, I2M, Marseille, France\\
  Sapienza Università di Roma, Dipartimento di Matematica, Roma, Italy\\
  CEREMADE, CNRS, Université Paris-Dauphine, Université PSL,
75016 Paris, France
  \BEMAIL{forien@ceremade.dauphine.fr},
  \url{https://www.ceremade.dauphine.fr/\string~forien}}
  \orcid{0000-0002-3146-644X}
  }


\AMSSUBJ{82B20}
\AMSSUBJSECONDARY{82B27 ; 60K35}

\SUBMITTED{December 21, 2021}
\ACCEPTED{January 5, 2024}

\ARXIVID{2110.07949} 
\HALID{hal-03382821} 

\VOLUME{0}
\YEAR{2023}
\PAPERNUM{0}
\DOI{10.1214/YY-TN}

\ABSTRACT{In~\cite{CG}, Cerf and Gorny constructed a model of self-organized criticality, by introducing an automatic control of the temperature parameter in the generalized Ising Curie-Weiss model.
The fluctuations of the magnetization of this spin model are of order~\smash{$n^{3/4}$} with a limiting law of the form~\smash{$C\exp(-x^4)$}, as in the critical regime of the Curie-Weiss model.

In this article, we build upon this model by replacing the mean-field interaction with a one-dimensional interaction with a certain range~$r_n$ which varies as a function of the number~$n$ of particles.
In the Gaussian case, we show that the self-critical behaviour observed in the mean-field case extends to interaction ranges~\smash{$r_n\gg n^{3/4}$} and we show that this threshold is sharp, with different fluctuations when the interaction range is of order of~$n^{3/4}$ or smaller than~$n^{3/4}$.}

\newcommand{\R}{\mathbb{R}}
\newcommand{\C}{\mathbb{C}}
\newcommand{\Z}{\mathbb{Z}}

\newcommand{\abs}[1]{\left| #1 \right|}
\newcommand{\acc}[1]{\left\{ #1 \right\}}
\newcommand{\Ent}[1]{\left\lfloor #1\right\rfloor}

\newcommand{\Sumj}{\sum_{j=1}^n}
\newcommand{\Sumju}{\sum_{j=1}^{n-1}}

\newcommand{\Prodju}{\prod_{j=1}^{n-1}}
\newcommand{\Prodi}{\prod_{i=1}^n}
\newcommand{\demi}{\frac{1}{2}}
\newcommand{\uc}{u^\star}
\newcommand{\limn}{\lim\limits_{n\to\infty}}
\newcommand{\cvloiinfty}{\underset{n \to +\infty}{\overset{\mathcal{L}}{\longrightarrow}}}
\newcommand{\eqninfty}{\stackrel{n\to\infty}{\sim}}
\newcommand{\alp}{\alpha_{n,j}}
\newcommand{\bet}{\beta_{n,j}}
\newcommand{\qquadet}{\qquad\text{and}\qquad}

\newcommand\numberthis{\addtocounter{equation}{1}\tag{\theequation}}

\newcommand{\guillemets}[1]{``#1''}

\makeatletter
\newcommand*{\biggg}[1]{{\hbox{$\left#1\vbox to20.5\p@{}\right.\n@space$}}}
\makeatother

\begin{document}



\section{Introduction}

\subsection{Definition of the model}

This article is devoted to the study of a one-dimensional spin model with long range interactions and with a self-adjusted temperature.
More precisely, we study a chain of~$n$ spins, with periodic boundary conditions, where each spin interacts with its~$2r_n$ nearest neighbours, with~\smash{$0<2r_n<n$}.
Thus, we consider the Hamiltonian
\begin{equation}
\label{eq:defHn}
H_n\ :\ (x_1,\,\ldots,\,x_n)\in\R^n\ \longmapsto\ -\frac{1}{2r_n}\sum_{i=1}^n\sum_{j=1}^{r_n} x_i x_{i+j}\,,
\end{equation}
where we use the convention~$x_{n+k}=x_k$ for all~$k\in\acc{1,\,\ldots,\,n}$.
Besides this Hamiltonian, the spins interact together through the following observable, which we think of as a \guillemets{self-adjusted temperature parameter} (which is not really a temperature parameter because it is a function of the spins):
\[T_n\ :\ (x_1,\,\ldots,\,x_n)\in\R^n\ \longmapsto\ \frac{x_1^2+\cdots+x_n^2}{n}\,.\]
We then study the following probability distribution on~$\R^n$:
\begin{equation}
\label{eq:defModele}
d\mu_n(x_1,\,\ldots,\,x_n)\ =\ 
\frac{1}{Z_n}\exp\left(-\frac{H_n(x_1,\,\ldots,\,x_n)}{T_n(x_1,\,\ldots,\,x_n)}\right)\mathbb{1}_{\acc{T_n>\,0}}\Prodi d\nu(x_i)\,,
\end{equation}
where~$\nu$ is the standard normal distribution, and~$Z_n$ is the normalization constant.
Note that the spins in our model are real-valued (if they were in~$\{-1,1\}$ then the temperature~$T_n$ would be constant).
In what follows, we are interested in the behaviour, when~$n$ tends to infinity, of the \guillemets{magnetization}
\[S_n(x_1,\,\ldots,\,x_n)\ =\ x_1+\cdots+x_n\,,\]
when~$(x_1,\,\ldots,\,x_n)$ is distributed according to~$\mu_n$.

\subsection{Results}
\label{sec:results}

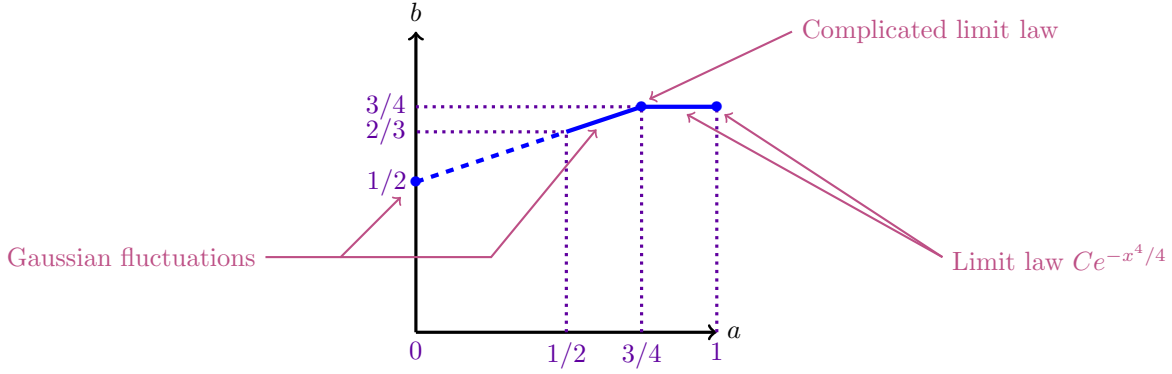
\begin{figure}
\centering
\begin{tikzpicture}
\draw[->, very thick] (0,0) -- (4,0) node[right]{$a$};
\draw[->, very thick] (0,0) -- (0,4) node[above]{$b$};
\draw[dotted, purple!50!blue, very thick] (2,0) node[below]{$1/2$} -- (2,2.6667) -- (0,2.6667) node[left]{$2/3$};
\draw[dotted, purple!50!blue, very thick] (3,0) node[below]{$3/4$} -- (3,3) -- (0,3) node[left]{$3/4$};
\draw[dotted, purple!50!blue, very thick] (4,0) node[below]{$1$} -- (4,3);
\draw[purple!50!blue] (0,2) node[left]{$1/2$};
\draw[purple!50!blue] (0,0) node[below]{$0$};
\draw[blue, ultra thick] (2,2.6667) -- (3,3) node{$\bullet$} -- (4,3) node{$\bullet$};
\draw[blue, ultra thick, dashed] (0,2) node{$\bullet$} -- (2,2.6667);
\draw[magenta!75!black, thick, ->] (-1.5,1) node[left]{Gaussian fluctuations} -- (-1,1) -- (-.2,1.8);
\draw[magenta!75!black, thick, ->] (-1,1) -- (1,1) -- (2.4,2.7);
\draw[magenta!75!black, thick, ->] (5.5,2) node[right]{Limit law $Ce^{-x^4/4}$} -- (3.6,2.9);
\draw[magenta!75!black, thick, ->] (5.5,2) -- (4.1,2.9);
\draw[magenta!75!black, thick, ->] (4.5,4) node[right]{Complicated limit law} -- (3.1,3.1);
\end{tikzpicture}
\caption{In the different regimes studied, if the interaction range is
of order~$r_n\sim n^a$, then the fluctuations of the
magnetization~$S_n$ are of order~$n^b$,
with~$b=(1/2+a/3)\wedge(3/4)$.}
\label{fig:results}
\end{figure}

We now describe the asymptotic behaviour of the magnetization in several regimes of the interaction range~$r_n$.
The results are summarized on figure~\ref{fig:results}.

\subsubsection{Self-critical behaviour in the long range case}

The following result indicates the behaviour of the magnetization in the long range case:

\begin{theorem}[Long range case]
\label{th:long_range}
If the interaction range~$r_n$ is such that~$r_n/n^{3/4}\to\infty$,
then under the law~$\mu_n$ defined by~(\ref{eq:defModele}), we have the convergence in distribution
\[\frac{S_n}{n^{3/4}}\ \cvloiinfty\ \frac{\sqrt{2}}{\Gamma(1/4)}\exp\left(-\frac{x^4}{4}\right)\,dx\,.\]
\end{theorem}

This interesting behaviour was already observed by Cerf and Gorny~\cite{GornyGaussien,CG} in the mean-field case, which roughly corresponds to an interaction range~$r_n=\Ent{n/2}$.
These fluctuations of the order of~$n^{3/4}$ and this limiting law correspond to the behaviour of the critical Ising-Curie-Weiss model~\cite{SG73,EN78,ENR80}.
Thus, the introduction of this self-adjusted temperature pushes the system to exhibit a critical-like behaviour, without the need to adjust the parameters of the system to precise values: this is why we talk of self-organized criticality (see paragraph~\ref{sec:motivation}).

\subsubsection{Threshold phenomenon}

A natural question is now: is this exponent~$3/4$ for the interaction range optimal?
The following Theorem provides a positive answer, showing that the limiting distribution changes when the interaction range~$r_n$ is of order~$n^{3/4}$:
\begin{theorem}[Threshold phenomenon]
\label{th:threshold}
If the interaction range~$r_n$ is such that
\[r_n\ \eqninfty\ \lambda n^{3/4}
\qquad\text{with}\qquad
\lambda>0\,,\]
if~$(Y_j)_{j\in\Z}$ is a family of i.i.d.\ standard normal variables, and if~$f$ is the density of the random variable
\begin{equation}
\label{eq:def_Z}
Z_\lambda\ =\ \sqrt{2}Y_0\,-\,\frac{3}{2\lambda^2\pi^2}\sum_{j\in\Z\setminus\acc{0}}\frac{Y_j^2}{j^2}\,,
\end{equation}
then under the law~$\mu_n$ defined by~(\ref{eq:defModele}), we have the convergence in distribution
\[\frac{S_n}{n^{3/4}}\ \cvloiinfty\ \frac{f\big(x^2\big)\,dx}{\int_{\R}f\big(t^2\big)\,dt}\,.\]
\end{theorem}
This Theorem indicates that there is a kind of phase transition phenomenon with respect to the interaction range.
The obtained limiting law results from the competition between the two terms in~(\ref{eq:def_Z}). When~$\lambda=\infty$, the second term disappears, leaving the distribution obtained in Theorem~\ref{th:long_range}.

\subsubsection{The finite range case}

When the interaction range is constant, the behaviour is very different from the long range case.
Indeed, the following Theorem shows that the phenomenon observed in the mean-field and long range cases does not occur in the finite range case, where the fluctuations of~$S_n$ are Gaussian:

\begin{theorem}[Finite range case]
\label{th:finite_range}
For every~$r\geqslant 1$, under the law~$\mu_n$ defined by~(\ref{eq:defModele}) with~$r_n=r$, we have the convergence in distribution
\[\frac{S_n}{\sqrt{n}}\ \cvloiinfty\ \mathcal{N}\big(0,\,\sigma_r^2\big)\,,\]
where~$\sigma_r>0$ is characterized by the equation
\begin{equation}
\label{eq:charac_sigma}
\int_0^1
dt\,\left(\frac{1}{\sigma_r^2}+1-\frac{1}{r}\sum_{m=1}^r \cos(2\pi mt)\right)^{-1}
\ =\ 1\,.
\end{equation}
In the nearest neighbour case~$r=1$, we have the explicit variance~$\sigma_1^2=\sqrt{2}+1$.
\end{theorem}

\subsubsection{An intermediate regime}

\label{sec:thm_intermediate}

The following theorem shows that the fluctuations of~$S_n$ become smaller than~$n^{3/4}$ when the interaction range becomes smaller than~$n^{3/4}$:

\begin{theorem}[Intermediate regime]
\label{th:intermediate}
If the interaction range satisfies~$r_n/\sqrt{n}\to\infty$ and~$r_n/n^{3/4}\to 0$ then, under the law~$\mu_n$ defined by~(\ref{eq:defModele}), we have the convergence in distribution
\[\frac{S_n}{{r_n}^{1/3}\sqrt{n}}\ \cvloiinfty\ \mathcal{N}\Bigg(0,\,\sqrt[3]{\frac{2}{3}}\Bigg)\,.\]
\end{theorem}

We have not been able to carry out the computations in the remaining
case~$r_n/\sqrt{n}\to 0$ but we expect that the variance~$\sigma_r^2$
defined by~(\ref{eq:charac_sigma}) satisfies~\smash{$\sigma_r\sim
(2/3)^{1/6}\,r^{1/3}$} when one lets~\smash{$r\to\infty$}, which leads us to conjecture
that Theorem~\ref{th:intermediate} remains true when~$r_n\to\infty$
and~$r_n/\sqrt{n}\to 0$. 
This guess corresponds to the dashed segment on
figure~\ref{fig:results}.

\subsection{Motivation: a Curie-Weiss model of self-organized criticality}
\label{sec:motivation}

The motivation which leads us to consider this model comes from the work of Cerf and Gorny~\cite{GornyGaussien,CG}, who studied a simple mean-field model of self-organized criticality which is constructed as a variant of the generalized Ising Curie-Weiss model.

The concept of self-organized criticality was coined in by the physicists Bak, Tang and Wiesenfeld in their seminal article~\cite{BTW87}, to explain the widespread presence of fractal structures in nature.
They observed that some physical systems present a \guillemets{critical-like} behaviour, with fractal structures and power-law correlations, without the need to finely tune a parameter (e.g., the temperature) to a critical value.
They called this phenomenon \guillemets{self-organized criticality}.
The important difference with ordinary phase transitions is that the critical regime, instead of being a very specific regime which only appears for a very precise value of the parameters of the system, becomes an attracting point, the system being \guillemets{forced} to look critical.
Several mathematical models of self-organized criticality have been studied, but these models are often quite complex and not easily tractable~\cite{These}.
For a broader review of the concept of self-organized criticality and of the controversies around it, we refer the reader to~\cite{Bak96, Frigg03, Pruessner12soc, WPCCJ16soc}.

To construct a simple toy model of self-organized criticality, Gorny started from the generalized Ising Curie-Weiss model and, following an idea explained by Sornette~\cite{SornettePerco}, he replaced the temperature parameter with a function of the spins, in order to introduce a kind of feedback from the configuration onto the temperature parameter.
Starting from a model of the form~$\exp(-H_n(\sigma)/T)$, with a phase transition for a critical temperature~$T_c$, the technique consists of replacing this temperature parameter~$T$ with a function~$T_n(\sigma)$, which tends to concentrate around the critical value~$T_c$ when~$n\to\infty$.

In~\cite{GornyGaussien}, Gorny constructed a mean-field model which almost corresponds to the case~$r_n=n-1$ of our model\footnote{In fact there is a small difference, with a factor~$1-1/n$ between our Hamiltonian and the one studied by Gorny, but this does not change the behaviour of the model.}, and he proved Theorem~\ref{th:long_range} in this setting, that is to say, he showed that the fluctuations of the magnetization are of order~$n^{3/4}$, with a limiting distribution of the form~$\exp(-x^4/4)$.
This corresponds to the behaviour observed in the critical regime of the Ising-Curie-Weiss model.

\subsection{Perspectives}

\subsubsection{More general distributions}

Our results only deal with the case of variables initially Gaussian, that is to say, we restrict ourselves to the standard normal distribution~$\nu$.
This restriction enables us to perform exact computations, as was done by Gorny in~\cite{GornyGaussien}, where he dealt with the Gaussian case of the model.
The work of Gorny was then extended to more general distributions of the spins in~\cite{CG} and~\cite{GV15}, still in a mean-field setting.

In fact, the model studied in the present article with an interaction range could also be defined with more general spin distributions.
But the Gaussian case with a varying interaction range turns out to be already challenging, and we can already observe an interesting threshold phenomenon in the interaction range.
Yet, we expect that our results should remain valid for a more general class of spin distributions.

To extend our results to more general distributions, it might be useful to use a Hubbard-Stratonovich transformation, as Gorny and Varadhan did in~\cite{GV15}, or to use the more general Stein's method~\cite{EL10}, as used in~\cite{DM20}.
In fact, our method uses a trick which consists in considering the self-normalized magnetization~\smash{$S_n/\sqrt{T_n}$} instead of~$S_n$, an idea which was already used in~\cite{GV15} and helped extending the result to more general distributions.
In our case, in addition to this, we benefit from the \guillemets{magical} fact that, in the Gaussian case, the self-normalized magnetization turns out to be independent from the temperature, which simplifies the computations.
Note that this fact remains true as long as the measure is spherically
symmetric: thus, our method easily extends to any spherically
symmetric base measure on~$\R^n$, as long as it is such that~$T_n\to
1$ in probability.
But this symmetry ingredient is not necessary, and the work of Gorny
and Varadhan shows that the method can be used in a more general setting.

\subsubsection{The intermediate regime}

As was explained in paragraph~\ref{sec:thm_intermediate}, we leave apart the regime between finite interaction range and~$r_n$ of order~$\sqrt{n}$, although we do not expect a different behaviour in this regime.
It would be interesting to check if something unexpected happens, or
if our limitation to~$\sqrt{n}$ is only an artefact of our method.
If this is the case, there might exist a smarter way to perform the computations which would allow to remove this limitation and to deal at the same time with the whole regime~\smash{$1\ll r_n\ll n^{3/4}$}.

\subsubsection{Different kinds of long range interactions}

The first motivation to study this model was to try to extend the construction of Cerf and Gorny, which was in a mean-field setting, to define a more geometrical model.
Thus, we studied the behaviour of the model in one of the simplest geometrical settings, namely a one-dimensional nearest neighbour interaction with periodic boundary conditions.
But, as shown by Theorem~\ref{th:finite_range}, it turns out that this model does not present the same critical-like asymptotic behaviour as observed for the mean-field model.
At this point, a natural question arises: if the interesting behaviour (with a limiting distribution of the form~$\exp(-x^4/4)$ and fluctuations of order~$n^{3/4}$) is observed in the mean-field case, but not when each spin only interacts with its two nearest neighbours, then what about the intermediate cases between these two extreme situations?

This naturally leads to consider a model with an intermediate interaction range.
But there are many different ways to interpolate between a mean-field interaction and a nearest neighbour interaction.
Generally speaking, one can consider a Hamiltonian~$H_n$ of the form
\[H_n(x_1,\,\ldots,\,x_n)\ =\ -\sum_{1\leqslant i,j\leqslant n}J(i,\,j)\, x_i\, x_j\,,\]
where the coupling constants~$J(i,\,j)=J(|i-j|)$ are decreasing functions of the distance separating the particles~$i$ and~$j$ (with periodic boundary conditions).
The behaviour of the model then depends on the decay rate of this coupling function.
In~\cite{ACCN88}, the key example of a coupling proportional to~$|i-j|^{-2}$ is studied, and it is proved that the resulting Ising model presents a phase transition, as well as models constructed with a slower-decaying coupling function.
Otherwise, if~$J(i,\,j)=o(|i-j|^{-2})$, then the obtained Ising model does not any more present a phase transition, which shows that a coupling of order~$|i-j|^{-2}$ plays a pivotal role for the appearance of a phase transition.

Another variant, called Kac-Ising models, consists of an interaction function (either with finite range or decreasing with the distance between spins) scaled by a factor~$\gamma$, and one studies the limit of these models when~$\gamma\to 0$~\cite{KUH63, BZ97}.

\subsubsection{Random couplings}

Another way to design intermediate models consists in drawing random couplings~$J(i,\,j)$.
In~\cite{BG93}, Bovier and Gayrard constructed such a model by taking for the~$J(i,\,j)$ independent Bernoulli variables of parameter~$p$, which amounts to considering the Ising model on an Erd{\H o}s-Rényi random graph.
This model exhibits different regimes characterized by different fluctuations of the sum of spins, depending on how the parameter~$p$ varies with the number~$n$ of particles.
These different regimes were studied by Kabluchko, L\"owe and Schubert~\cite{KLS19a, KLS19b} who proved in particular that, for a critical temperature and a parameter~$p_n$ chosen such that~$p_n/n^{-3/4}\to\infty$, the behaviour resembles that of the critical Ising Curie-Weiss model, i.e., the sum of the spins is of order~$n^{3/4}$ with fluctuations of the form~$C\exp(-\lambda x^4)$.
When the parameter~$p_n$ becomes of order~$n^{-3/4}$, still at the critical temperature, the limiting distribution changes, and a quadratic term appears besides the term in~$x^4$.
If~$p_n=o(n^{-3/4})$ then this quadratic term dominates, which results in Gaussian fluctuations of the sum of spins.

This approach was generalized by Deb and Mukherjee~\cite{DM20}, who studied the fluctuations of an Ising model defined on a more general set of graphs.
Under certain conditions of homogeneity and connectivity, they obtain the same fluctuations as in the mean-field model, when the mean degree~$r_n$ in the graph satisfies~$r_n/(n\ln n)^{1/3}\to\infty$ in the supercritical regime,~$r_n/\sqrt{n}\to\infty$ in the subcritical regime, or~$r_n/\sqrt{n}\ln n\to\infty$ at the critical point.

In our case, we choose an interaction of each spin with its~$2r_n$ nearest neighbours, where~$r_n$ is a parameter which evolves with~$n$.
This corresponds to a coupling function of the form
\[J(i,\,j)
\ =\ \frac{1}{4r_n}\times\left\{\begin{aligned}
&1\quad\text{ if}\ j\in\acc{i-r_n,\,\ldots,\,i-1}\cup\acc{i+1,\,\ldots,\,i+r_n}+n\Z\,,\\
&0\quad\text{ otherwise.}
\end{aligned}\right.\]

Let us recall that our model is different from the aforementioned models, because the spins are not valued in~$\acc{-1,\,+1}$ but are real-valued, and because our self-adjusted temperature in fact induces an interaction between all the spins.
Therefore, it is not \textit{a priori} evident which scale of~$r_n$ is relevant to observe a change of behaviour.

\subsection{Sketch of the proof and organization of the paper}

\subsubsection{Diagonalization of the interaction Hamiltonian}

Our starting point to study the model is to diagonalize the interaction Hamiltonian.
The matrix of the quadratic form~$H_n$ given by~(\ref{eq:defHn}) is a symmetric circulant matrix, which writes
\[\mathfrak{M}\big(H_n\big)
\ =\ -\frac{1}{4r_n}\sum_{m=1}^{r_n}\big(J^m+J^{-m}\big)\,,\]
where the matrix~$J$ is given by
\[J\ =\ \begin{pmatrix}\,
0 & 1 & 0 & \cdots & 0\\
0 & 0 & 1 & \ddots & \vdots\\
\vdots & & \ddots & \ddots & 0\\
0 & & & 0 & 1\\
1 & 0 & \cdots & \cdots & 0
\,\end{pmatrix}
\ =\ \Big(\mathbb{1}_{\{j=i+1\text{ mod.\ }n\}}\Big)_{1\leqslant i,j\leqslant n}\,.\]
The matrix~$J$ is diagonalizable, with eigenvalues
\[\mathrm{sp}(J)\ =\ \Big\{\,e^{2ij\pi/n}\ :\ j\in\acc{1,\,\ldots,\,n}\,\Big\}\,.\]
Therefore, for every~$m\geqslant 1$, we have
\[\mathrm{sp}\big(J^m+J^{-m}\big)\ =\ \acc{\,2\cos\left(\frac{2jm\pi}{n}\right)\ :\ j\in\acc{1,\,\ldots,\,n}\,}\,,\]
and thus the eigenvalues of our Hamiltonian~$H_n$ are
\[\mathrm{sp}(H_n)\ =\ \bigg\{\,-\frac{\alpha_{n,1}}{2},\,\ldots,\,-\frac{\alpha_{n,n}}{2}\,\bigg\}\,,\]
where
\begin{equation}
\label{eq:alpha_cosine}
\forall j\in\acc{1,\,\ldots,\,n}\qquad
\alp\ =\ \frac{1}{r_n}\sum_{m=1}^{r_n}{\cos\left(\frac{2jm\pi}{n}\right)}\,.
\end{equation}
Also, for every~$j\in\{1,\,\ldots,\,n\}$, we define~$\bet=1-\alp$.
In section~\ref{sec:eigenvalues}, we gather a series of estimates on these eigenvalues that are used throughout the article.

\subsubsection{Change of variables}

Let~$P\in\mathcal{O}_n(\R)$ be an orthonormal matrix such that the
matrix of the quadratic form~$H_n$ writes~$P^{-1}DP$, where~$D$ is the
diagonal matrix with
coefficients~$-\alpha_{n,1}/2,\,\ldots,\,-\alpha_{n,n}/2$ on the
diagonal.
Performing the change of variables
\[\varphi\ :\ (x_1,\,\ldots,\,x_n)\in\R^n\ \longmapsto\ 
(y_1,\,\ldots,\,y_n)
\ =\ \left(\,\sum_{k=1}^n{P_{j,k}x_k}\,\right)_{1\leqslant j\leqslant n}\,,\]
the Hamiltonian~$H_n$ and the self-adjusted temperature~$T_n$ become
\[H_n(x_1,\,\ldots,\,x_n)\ =\ -\demi\Sumj{\alp y_j^2}\,,
\qquadet
T_n(x_1,\,\ldots,\,x_n)\ =\ \frac{1}{n}\Sumj{y_j^2}\,.\]
To see what happens to the sum~$S_n=x_1+\cdots+x_n$ of the spins, note that~$\alpha_{n,n}=1$ and that this eigenvalue~$-\alpha_{n,n}/2=-1/2$ of~$H_n$ corresponds to the eigenvector~$(1/\sqrt{n},\,\ldots,\,1/\sqrt{n})$, whence~\smash{$P_{n,k}=1/\sqrt{n}$} for all~$k\in\acc{1,\,\ldots,\,n}$.
Therefore, we have
\[(y_1,\,\ldots,\,y_n)\ =\ \varphi(x_1,\,\ldots,\,x_n)
\qquad\Rightarrow\qquad
S_n\ =\ x_1+\cdots+x_n\ =\ \sqrt{n}\,y_n\,.\]
Also, the change of variable being orthonormal, we have
\[\Prodi d\nu(x_i)
\ =\ d\nu^{\otimes n}\varphi^{-1}(y_1,\,\ldots,\,y_n)
\ =\ \Prodi d\nu(y_i)\,,\]
which allows us to forget the variables~$x_1,\,\ldots,\,x_n$ and to work only with the new variables~$y_1,\,\ldots,\,y_n$ to study the limiting behaviour of~$S_n$, which is now~$S_n=\sqrt{n}\,y_n$.

\subsubsection{A competition between two terms}
Using the fact that~$\alpha_{n,n}=1$, we may write
\begin{equation}
\label{eq:competition}
H_n(x_1,\,\ldots,\,x_n)\ =\ -\demi\Sumj\alp y_j^2
\ =\ -\frac{y_n^2}{2}-\demi\Sumju\alp y_j^2
\ =\ -\frac{S_n^2}{2n}-\demi\Sumju\alp y_j^2\,.
\end{equation}
The first term exactly corresponds to the Hamiltonian in the mean-field model of Cerf and Gorny, as mentioned in section~\ref{sec:motivation}.
In the mean-field case~$2r_n=n-1$ we have~\smash{$\alp=-1/(2r_n)$} for every~$j\neq n$, so that the second term almost disappears (in fact it contributes as~\smash{$-S_n^2/(4nr_n)$} but this has no significant effect).
Hence, in this case the model behaves as if there were only the first term~$-S_n^2/(2n)$.

For a general interaction range~$r_n$, there is a competition between the two terms in~(\ref{eq:competition}).
As long as~\smash{$r_n/n^{3/4}\to\infty$}, we show that the first term prevails, whence a behaviour similar to the mean-field case.
The interesting change of behaviour intervenes when the second term becomes big enough to disturb the first term, and it turns out that this happens when~$r_n$ is of order~\smash{$n^{3/4}$}.

\subsubsection{Independence of the temperature}
The following important observation simplifies significantly the study of the model.

\begin{lemma}[Independence of the temperature]
\label{lm:indep_T}
Let~$n\geqslant 1$.
If~$Y_1,\,\ldots,\,Y_n$ are i.i.d.\ standard normal variables (we
denote by~$\nu^{\otimes n}$ their joint distribution), then the
variable~\smash{$T_n=(Y_1^2+\cdots+Y_n^2)/n$} is independent of~\smash{$\big(Y_1/\sqrt{T_n},\,\ldots,\,Y_n/\sqrt{T_n}\big)$}.
For every measurable and bounded function~$g:\R^n\to\R$, we have
\begin{multline*}
\nu^{\otimes n}\Bigg[\,g\left(\frac{Y_1}{\sqrt{T_n}},\,\ldots,\,\frac{Y_n}{\sqrt{T_n}}\right)\,\Bigg]
\ =\ \frac{\Gamma(n/2)}
{\pi^{n/2} n^{(n-2)/2}}\\
\times
\int_{\R^{n-1}}
\frac{dz_1\,\ldots\,dz_{n-1}\,\mathbb{1}_{\{z_1^2+\cdots+z_{n-1}^2<n\}}}
{\sqrt{n-z_1^2-\cdots-z_{n-1}^2}}\,
g_s\bigg(z_1,\,\ldots,\,z_{j-1},\,\sqrt{n-z_1^2-\cdots-z_{n-1}^2}\bigg)\,,
\end{multline*}
where~$g_s$ is the even part of~$g$ with respect to the last variable, that is to say,
\begin{equation}
\label{eq:defGs}
g_s\,:\,
(z_1,\,\ldots,\,z_n)\in\R^n
\ \longmapsto\ \frac{g(z_1,\,\ldots,\,z_{n-1},\,z_n)+g(z_1,\,\ldots,\,z_{n-1},\,-z_n)}{2}\,.
\end{equation}
\end{lemma}

This Lemma is proved in section~\ref{sec:indep_temperature}.
It follows from this result that, under~$\nu^{\otimes n}$, that is to say when the variables~$Y_j$ are i.i.d.\ standard Gaussian variables, the variable~$T_n$ is independent of
\[\exp\left(-\frac{H_n}{T_n}\right)
\ =\ \exp\Bigg(\,\demi\Sumj\alp\frac{Y_j^2}{T_n}\,\Bigg)\,.\]
Therefore, for every bounded and measurable function~$f:\R\to\R$, we have
\begin{align*}
\mu_n\Big[f\big(T_n\big)\Big]
\ &=\ \frac{1}{Z_n}\nu^{\otimes n}\left[\,\exp\left(-\frac{H_n}{T_n}\right)f\big(T_n\big)\,\right]\\
\ &=\ \frac{1}{Z_n}\nu^{\otimes n}\left[\,\exp\left(-\frac{H_n}{T_n}\right)\,\right]\nu^{\otimes n}\Big[f\big(T_n\big)\Big]
\ =\ \nu^{\otimes n}\Big[f\big(T_n\big)\Big]\,.
\end{align*}
Thus, the variable~$T_n$ has the same distribution in our model~$\mu_n$ as if the variables~$Y_j$ were i.i.d.\ standard normal variables, so the law of large numbers implies the following fact:

\begin{corollary}[Behaviour of the temperature]
\label{cor:temperature}
Under~$\mu_n$, we have the convergence in probability~$T_n\to 1$.
\end{corollary}

\subsubsection{Slutsky's Lemma and self-normalized magnetization}

Now that we know that~$T_n\to 1$ in probability, by virtue of Slutsky's Lemma (Theorem~3.9 in~\cite{Billingsley99}), to prove the convergences in distribution announced in section~\ref{sec:results}, it is enough to show the same results for the self-normalized magnetization~$S_n/\sqrt{T_n}$ instead of~$S_n$, as was done in~\cite{GV15}.
Thus, there now remains to study the behaviour of this variable~$S_n/\sqrt{T_n}$.
Using the formula for the distribution of~\smash{$\big(Y_1/\sqrt{T_n},\,\ldots,\,Y_n/\sqrt{T_n}\big)$} given by Lemma~\ref{lm:indep_T} and plugging it into the definition of our model, we easily obtain:

\begin{lemma}[Distribution of the self-normalized magnetization]
\label{lm:distrib_magn}
For every bounded and measurable function~$g:\R\to\R$ and every~$n\geqslant 3$, we have
\[\mu_n\Bigg[g\left(\frac{S_n}{\sqrt{T_n}}\right)\Bigg]
\ =\ C_n\,
\mathbb{E}\bigg[\widetilde{g}\Big(n\big(n-A_n\big)\Big)\mathbb{1}_{\{A_n<n\}}\bigg]\]
with~$A_n=Z_1^2+\cdots+Z_{n-1}^2$, where the variables~$Z_j$ are independent with~\smash{$Z_j\sim\mathcal{N}\big(0,\,1/\beta_{n,j}\big)$}, and where the constant~$C_n$ and the function~$\widetilde{g}$ are given by
\[C_n
\ =\ 
\frac{\Gamma(n/2)e^{(n-3)/2}2^{n/2}}
{n^{n/2}Z_n}\,
\Prodju\frac{1}{\sqrt{\bet}}
\qquadet
\widetilde{g}
\,:\,x\in(0,\,\infty)\ \longmapsto\ 
\frac{g\big(\sqrt{x}\big)+g\big(-\sqrt{x}\big)}
{2\sqrt{x}}\,.\]
\end{lemma}

The above Lemma has the following corollary, which will allow us to deduce the asymptotic behaviour of~$S_n/\sqrt{T_n}$ from the asymptotic behaviour of the auxiliary variable~$n-A_n$:

\begin{corollary}[Relating the self-normalized magnetization to the
behaviour of~$A_n$]
\label{cor:relating_behaviours}
If, for a certain sequence~$(v_n)_{n\geqslant 3}$ of positive real numbers, the random variable~$n(n-A_n)/v_n^2$, conditioned to be positive,
converges in distribution to a random variable with density~$f$ with respect to the Lebesgue measure on~$\R$, then we have the convergence in distribution
\[\frac{S_n}{v_n\sqrt{T_n}}
\ \cvloiinfty\ 
\frac{f\big(x^2\big)\,dx}{\int_{\R}f\big(t^2\big)\,dt}\,.\]
\end{corollary}

The proof of Lemma~\ref{lm:distrib_magn} and of this Corollary~\ref{cor:relating_behaviours} are given in section~\ref{sec:relating_behaviours}.

\subsubsection{Behaviour of the auxiliary variable in the various regimes}

Now, there only remains to study the behaviour of the random variable~$n-A_n$ conditioned to be positive, in the various regimes of the interaction range~$r_n$, and with a suitable scaling~$v_n$.
This behaviour is described by the following four Lemmas.
In what follows, we keep the notation of Lemma~\ref{lm:distrib_magn} above, writing~\smash{$A_n=Z_1^2+\cdots+Z_{n-1}^2$}, where the variables~$Z_j$ are independent with~\smash{$Z_j\sim\mathcal{N}\big(0,\,1/\beta_{n,j}\big)$}, recalling that~$\bet=1-\alp$.

\begin{lemma}[Behaviour of~$A_n$ in the long range case]
\label{lm:long_range}
If~$r_n/n^{3/4}\to\infty$, then we have the convergence in distribution
\[\frac{n-A_n}{\sqrt{n}}\ \cvloiinfty\ \mathcal{N}(0,\,2)\,.\]
\end{lemma}

\begin{lemma}[Behaviour of~$A_n$ at the threshold]
\label{lm:threshold}
If~$r_n\sim\lambda n^{3/4}$ with~$\lambda>0$, then we have the convergence in distribution
\[\frac{n-A_n}{\sqrt{n}}\ \cvloiinfty\ Z_\lambda\,,\]
where~$Z_\lambda$ is the variable defined by~(\ref{eq:def_Z}).
\end{lemma}

\begin{lemma}[Behaviour of~$A_n$ in the finite range case]
\label{lm:finite_range}
If~$r_n=r$ is constant, then, conditioned on the event~$\{A_n<n\}$, the variable~$n-A_n$
converges in distribution to an exponential variable with parameter~$1/(2\sigma_r^2)$, where~$\sigma_r$ is characterized by~(\ref{eq:charac_sigma}).
\end{lemma}

\begin{lemma}[Behaviour of~$A_n$ in the intermediate regime]
\label{lm:intermediate}
If the range~$r_n$ is such that~$r_n/\sqrt{n}\to\infty$ and~$r_n/n^{3/4}\to 0$ then, conditioned on the event~$\{A_n<n\}$, the random variable~\smash{$(n-A_n)/r_n^{2/3}$} converges in distribution to an exponential variable with parameter~$3^{1/3}/2^{4/3}$.
\end{lemma}

Together with Corollaries~\ref{cor:temperature} and~\ref{cor:relating_behaviours} and Slutsky's Lemma, the four above Lemmas easily imply our main results, namely Theorems~\ref{th:long_range},~\ref{th:threshold},~\ref{th:finite_range} and~\ref{th:intermediate}.

Note that, as explained above, since we expect
that~$r^{2/3}/(2\sigma_r^2)\to 3^{1/3}/2^{4/3}$ when one lets~$r\to\infty$, we
conjecture that the hypothesis~$r_n/\sqrt{n}\to\infty$ in
Lemma~\ref{lm:intermediate} is not optimal. It might be that this
hypothesis could be improved to~$r_n\to\infty$, but this would require
a more delicate analysis (see the proof of
Lemma~\ref{lm:intermediate} in section~\ref{sec:Intermediate} and the remark after
the proof).

To study the asymptotic behaviour of the auxiliary variable~$n-A_n$ in these various regimes, our common starting point is its characteristic function given by the following Lemma, which immediately follows from the formula for the characteristic function of the chi-squared distribution (see for example~\cite{Feller}).

\begin{lemma}[Fourier transform of~$A_n$]
\label{lm:fourier}
For every~$n\geqslant 1$, the characteristic function of~$n-A_n$ writes
\[\mathbb{E}\Big[e^{iu(n-A_n)}\Big]
\ =\ e^{iun}\Prodju \mathbb{E}\Big[e^{-iuZ_j^2}\Big]
\ =\ \exp\big[\varphi_n(u)\big]\,,\]
where~$\varphi_n$ is the function
\begin{equation}
\label{eq:def_phi}
\varphi_n
\,:\,u\in U_n\ \longmapsto\ 
iun-\demi\Sumju \ln \left(1+\frac{2iu}{\beta_{n,j}}\right)\,,
\end{equation}
which can be defined on the complex domain
\begin{equation}
\label{eq:def_un}
U_n
\ =\ \Big\{u\in\C\ :\ 2\,\mathfrak{Im}\,u\,<\,\mathrm{min}\big(\beta_{n,1},\,\ldots,\,\beta_{n,n-1}\big)\,\Big\}\,,
\end{equation}
using the following determination of the logarithm:
\[\ln\,:\,z=x+iy\in\C\setminus(-\infty,0]
\ \longmapsto\ \demi\ln(x^2+y^2)+2i\arctan\left(\frac{y}{x+\sqrt{x^2+y^2}}\right)\,.\]
\end{lemma}

To obtain Lemmas~\ref{lm:long_range} and~\ref{lm:threshold}, we simply compute the limit of the characteristic function with a suitable scaling.
This is done in sections~\ref{sec:LongRange} and~\ref{sec:Threshold}.
In these cases, the conditioning on the event~$\{A_n<n\}$ is completely transparent, so that we may study the convergence in distribution without this conditioning.

\subsubsection{Saddle-point method}

The method in the regimes when~\smash{$r_n=o\big(n^{3/4}\big)$} is less straightforward because, in this case, it turns out that~$\mathbb{P}(A_n<n)\to 0$ when~$n\to\infty$, that is to say, we have to study some large deviations of the variable~$n-A_n$.

To do so, we use the Fourier transform to obtain an exact integral formula for the density of~$n-A_n$ and, to understand the limiting behaviour of the integral, we use the saddle-point method.
This method (see for example~\cite{Copson}) consists in a judicious change of integration contour in the complex plane in order to obtain an integrand with an appropriate limiting behaviour.

We perform this change of contour in section~\ref{sec:MovingContour},
before studying the limiting behaviour of the obtained integral, first
in the case of~$r_n$ constant in section~\ref{sec:FiniteRange} and
then in the intermediate regime of Theorem~\ref{th:intermediate} in
section~\ref{sec:Intermediate}.

This idea to use the saddle-point method to study a long range Ising model, after having diagonalized the interaction matrix, was already presented by Canning in a series of publications~\cite{Canning1992class, Canning92saddle, Canning1993generalized}.
However, Canning only discusses the saddle-point method in cases where the interaction matrix has a finite rank, and he does not give rigorous bounds on the precision of the obtained approximation.
In our computations, it turns out that the proof of the convergence in distribution requires a precise control of the asymptotic behaviour of the eigenvalues~$\alp$.
The analytical results on these eigenvalues that we need are gathered in the next section.

\section{Preliminaries}

\subsection{Estimates on the eigenvalues}
\label{sec:eigenvalues}

We now state some preliminary estimates on the eigenvalues which are needed in the various regimes of the interaction range.
These eigenvalues~$(\alp)_{1\leqslant j\leqslant n}$ were defined by~(\ref{eq:alpha_cosine}).
Note that we have~$\alpha_{n,n}=1$ and, for every~$j\in\acc{1,\,\ldots,\,n-1}$,
\begin{align*}
\alp
\ =\ \frac{1}{r_n}\sum_{m=1}^{r_n}{\cos\left(\frac{2jm\pi}{n}\right)}
&\ =\ \frac{\cos\big(j(r_n+1)\pi/n\big)\sin\big(jr_n\pi/n\big)}{r_n\sin\big(j\pi/n\big)}\numberthis\label{eq:alpha_cos_sin_sin}\\
&\ =\ \frac{1}{2r_n}\left(\frac{\sin\big((2r_n+1)j\pi/n\big)}{\sin\big(j\pi/n\big)}-1\right)\,.\numberthis\label{eq:alpha_sin_sin}
\end{align*}
The eigenvalues are symmetric with respect to~$j$, that is to say~$\alp=\alpha_{n,n-j}$ for every~$j\in\acc{1,\,\ldots,\,n-1}$.
It follows from~(\ref{eq:alpha_cos_sin_sin}) that
\begin{equation}
\label{eq:upper_bound_alpha}
\forall j\in\acc{\,1,\,\ldots,\,\Ent{\frac{n}{2}}\,}\qquad
\abs{\alp}
\ \leqslant\ \frac{1}{r_n\sin(j\pi/n)}
\ \leqslant\ \frac{n}{2r_n j}\,.
\end{equation}
For every~$j\in\acc{1,\,\ldots,\,n}$, we write~$\bet=1-\alp$.
The proofs of the following technical Lemmas are deferred to the appendix~\ref{sec:study_eigenvalues}.

\begin{lemma}[Upper bounds on the eigenvalues]
\label{lm:bound_eigenvalues}
We have the upper bounds
\[\sum_{j=1}^{n}\abs{\alp}
\ =\ O\left(\frac{n\ln n}{r_n}\right)
\qquadet
\sum_{j=1+\Ent{n/r_n}}^{n-\Ent{n/r_n}-1}\abs{\frac{1}{\bet}-1}
\ =\ O\left(\frac{n\ln n}{r_n}\right)\,.\]
\end{lemma}

\begin{lemma}[Asymptotics of the eigenvalues]
\label{lm:asymptotics_eigenvalues}
If~$r_n\to\infty$, then there exists~$K>0$ such that, for~$n\geqslant 3$ and~$1\leqslant j\leqslant\Ent{n/2}$,
\[\abs{\frac{(2r_n)^2}{n^2\bet}
-\frac{6}{\pi^2 j^2}
}\ \leqslant\ \frac{K}{r_n j^2}
+\frac{K{r_n}^2}{n^2}\,.\]
\end{lemma}

\subsection{Independence of the temperature: proof of Lemma~\ref{lm:indep_T}}
\label{sec:indep_temperature}

\begin{proof}[Proof of Lemma~\ref{lm:indep_T}]
Let~$f:\R\to\R$ and~$g:\R^n\to\R$ be two measurable and bounded functions, and define, for every~$n\geqslant 1$,
\[I_n\ =\ \nu^{\otimes n}\Bigg[\,f(T_n)
\,g\left(\frac{Y_1}{\sqrt{T_n}},\,\ldots,\,\frac{Y_n}{\sqrt{T_n}}\right)\,\Bigg]\,.\]
Considering the function~$g_s$ defined as in~(\ref{eq:defGs}), we can write
\begin{align*}
I_n
\ &=\ 
\int_{\R^n}\frac{dy_1\,\ldots\,dy_n}{(2\pi)^{n/2}}\,
f\Bigg(\frac 1 n \Sumj y_j^2\Bigg)\,
g\left(\frac{y_1\sqrt{n}}{\sqrt{\Sumj y_j^2}},\,\ldots,\,
\frac{y_n\sqrt{n}}{\sqrt{\Sumj y_j^2}}\right)
\exp\Bigg(-\Sumj\frac{y_j^2}{2}\Bigg)
\\
\ &=\ 
\frac{2}{(2\pi)^{n/2}}
\int_{\R^{n-1}}dy_1\,\ldots\,dy_{n-1}
\int_0^\infty dy_n\,
f\Bigg(\frac 1 n \Sumj y_j^2\Bigg)\\
&\qquad\times
g_s\left(\frac{y_1\sqrt{n}}{\sqrt{\Sumj y_j^2}},\,\ldots,\,
\frac{y_n\sqrt{n}}{\sqrt{\Sumj y_j^2}}\right)
\,\exp\Bigg(-\Sumj\frac{y_j^2}{2}\Bigg)\,.
\end{align*}
We now perform the change of variable~$y_n=\sqrt{nt-y_1^2-\cdots -y_{n-1}^2}$, which leads to
\begin{multline*}
I_n
\ =\ \frac{2}{(2\pi)^{n/2}}
\int_{\R^{n-1}}dy_1\,\ldots\,dy_{n-1}
\int_0^\infty
dt\, 
\frac{n\,\mathbb{1}_{\{nt>y_1^2+\cdots+y_{n-1}^2\}}}{2\sqrt{nt-y_1^2-\cdots -y_{n-1}^2}}
f(t)
\,\\
\times\ g_s\left(\frac{y_1}{\sqrt{t}},\,\ldots,\,\frac{y_{n-1}}{\sqrt{t}}\,,
\sqrt{\frac{nt-y_1^2-\cdots-y_{n-1}^2}{t}}\right)
\exp\left(-\frac{nt}{2}\right)\,.
\end{multline*}
Swapping the summations using Fubini's Theorem and writing~$y_j=z_j\sqrt{t}$ for every~$j<n$, we obtain
\begin{multline*}
I_n
\ =\ \frac{n}{(2\pi)^{n/2}}
\left(\int_0^\infty dt\,t^{(n-2)/2}e^{-nt/2}f(t)\right)\\
\times\ \Bigg(\int_{\R^{n-1}}
\frac{dz_1\,\ldots\,dz_{n-1}\,\mathbb{1}_{\{z_1^2+\cdots+z_{n-1}^2<n\}}}
{\sqrt{n-z_1^2-\cdots-z_{n-1}^2}}\,
g_s\Big(z_1,\,\ldots,\,z_{n-1},\,
\sqrt{n-z_1^2-\cdots-z_{n-1}^2}\Big)\Bigg)\,,
\end{multline*}
which shows that the variables~$T_n$ and~\smash{$\big(Y_1/\sqrt{T_n},\,\ldots,\,Y_n/\sqrt{T_n}\big)$} are independent, and yields the formula for the distribution of the latter variable which appears in the statement of the Lemma.
\end{proof}

\subsection{Relating the self-normalized magnetization to the auxiliary variable}
\label{sec:relating_behaviours}

This section is devoted to the proof of Lemma~\ref{lm:distrib_magn} and of Corollary~\ref{cor:relating_behaviours}.

\begin{proof}[Proof of Lemma~\ref{lm:distrib_magn}]
Let~$g:\R\to\R$ be a measurable and bounded function, and let~$n\geqslant 3$.
Recalling the definition~(\ref{eq:defModele}) of our model and defining
\[g_s\,:\,x\in\R\ \longmapsto\ 
\frac{g(x)+g(-x)}{2}\,,\]
we can write
\begin{align*}
\mu_n\Bigg[g\left(\frac{S_n}{\sqrt{T_n}}\right)\Bigg]
&\ =\ \frac{1}{Z_n}\nu^{\otimes n}\Bigg[\,g\left(\frac{S_n}{\sqrt{T_n}}\right)\,
\exp\left(-\frac{H_n}{T_n}\right)\,\Bigg]
\\
&\ =\ \frac{1}{Z_n}\nu^{\otimes n}\Bigg[\,g\left(\frac{Y_n\sqrt{n}}{\sqrt{T_n}}\right)\,
\exp\Bigg(\demi\Sumj\alp\frac{Y_j^2}{T_n}\Bigg)\,\Bigg]
\\
&\ =\ \frac{\Gamma(n/2)}
{\pi^{n/2} n^{(n-2)/2}Z_n}
\int_{\R^{n-1}}
\frac{dz_1\,\ldots\,dz_{n-1}}
{\sqrt{n-z_1^2-\cdots-z_{n-1}^2}}\,
\mathbb{1}_{\{z_1^2+\cdots+z_{n-1}^2<n\}}
\\
&\hspace{-2em}\times\ g_s\left(\sqrt{n\big(n-z_1^2-\cdots-z_{n-1}^2\big)}\right)\,
\exp\Bigg(\demi\Sumju \alp z_j^2
+\frac{n-z_1^2-\cdots-z_{n-1}^2}{2}\Bigg)
\,,
\end{align*}
using Lemma~\ref{lm:indep_T}.
With~$\widetilde{g}$ and~$A_n$ defined as in the statement, this can
be rewritten
\begin{align*}
\mu_n\Bigg[g\left(\frac{S_n}{\sqrt{T_n}}\right)\Bigg]
&\ =\ 
\frac{\Gamma(n/2)e^{n/2}}
{\pi^{n/2} n^{(n-3)/2}Z_n}
\int_{\R^{n-1}}
dz_1\,\ldots\,dz_{n-1}\,
\mathbb{1}_{\{z_1^2+\cdots+z_{n-1}^2<n\}}
\\
&\qquad\qquad\times\ 
\widetilde{g}\Big(n\big(n-z_1^2-\cdots-z_{n-1}^2\big)\Big)\,
\exp\Bigg(-\demi\Sumju \bet z_j^2\Bigg)
\\
&\ =\ 
\frac{\Gamma(n/2)e^{n/2}2^{n/2}}
{n^{(n-3)/2}Z_n}\,
\Prodju\frac{1}{\sqrt{\bet}}\,
\mathbb{E}\bigg[\widetilde{g}\Big(n\big(n-A_n\big)\Big)\mathbb{1}_{\{A_n<n\}}\bigg]\,,
\end{align*}
which concludes the proof of the Lemma.
\end{proof}

We now briefly explain how Corollary~\ref{cor:relating_behaviours} follows from this Lemma:

\begin{proof}[Proof of Corollary~\ref{cor:relating_behaviours}]
Assume that~$n(n-A_n)/v_n^2$, conditioned to be positive,
converges in distribution to a random variable with density~$f$ with respect to the Lebesgue measure on~$\R$, and let~$g:\R\to\R$ be a bounded and continuous function.
Following Lemma~\ref{lm:distrib_magn}, for~$n\geqslant 3$, we have
\begin{align*}
\mu_n\Bigg[g\left(\frac{S_n}{v_n\sqrt{T_n}}\right)\Bigg]
&\ =\ 
\frac{C_n}{v_n}\,
\mathbb{E}\Bigg[\widetilde{g}\left(\frac{n}{v_n^2}\big(n-A_n\big)\right)
\mathbb{1}_{\{A_n<n\}}\Bigg]
\\
&\!\!\eqninfty\ 
\frac{C_n\mathbb{P}[A_n<n]}{v_n}\,
\int_0^\infty
dx\,f(x)\,\widetilde{g}(x)
\,.
\end{align*}
Then, notice that
\[
\int_0^\infty
dx\,f(x)\,\widetilde{g}(x)
\ =\ 
\int_0^\infty
dx\,f(x)\frac{g_s\big(\sqrt{x}\big)}{\sqrt{x}}
\ =\ 
2
\int_0^\infty
dy\,f\big(y^2\big)\,g_s(y)
\ =\ 
\int_\R
dy\,f\big(y^2\big)\,g(y)
\,.
\]
To conclude the proof, there only remains to check that
\[\limn\,
\frac{C_n\mathbb{P}[A_n<n]}{v_n}\,
\ =\ 
\left(
\int_\R
dy\,f\big(y^2\big)
\right)^{-1}
\,,\]
which follows by taking~$g$ to be the constant function~$g=1$.
\end{proof}

\section{Study of the asymptotic behaviour of the auxiliary variable}

We now prove Lemmas~\ref{lm:long_range},~\ref{lm:threshold},~\ref{lm:finite_range} and~\ref{lm:intermediate}, that is to say, we study the behaviour of the random variable~$n-A_n$ in the various regimes of the interaction range.

\subsection{The long range case: proof of Lemma~\ref{lm:long_range}}
\label{sec:LongRange}

To prove the announced convergence in law, we simply compute the characteristic function and check that it tends to that of the normal distribution.

\begin{proof}[Proof of Lemma~\ref{lm:long_range}]
Assume that~$r_n/n^{3/4}\to\infty$.
It follows from Lemma~\ref{lm:fourier} that the characteristic function of~$(n-A_n)/\sqrt{n}$ writes
\begin{equation}
\label{eq:characNorm}
\mathbb{E}\Bigg[\exp\left(iu\,\frac{n-A_n}{\sqrt{n}}\right)\Bigg]
\ =\ \exp\Bigg[\varphi_n
\left(\frac{u}{\sqrt{n}}\right)
\Bigg]\,,
\end{equation}
where~$\varphi_n$ is the function defined by~(\ref{eq:def_phi}).
It follows from Taylor's formula that, for every~$y\in\R$,
\[\abs{\,
\ln (1+iy)-iy-\frac{y^2}{2}\,}
\ \leqslant\ \frac{1}{6}\,\sup_{t\in[0,1]}\,
\frac{2\,|y|^3}{\abs{1+ity}^3}
\ \leqslant\ \frac{|y|^3}{3}\,.\]
Plugging this into~(\ref{eq:def_phi}), we deduce that, for every~$u\in\R$,
\begin{equation}
\label{eq:DLSu}
\varphi_n
\left(\frac{u}{\sqrt{n}}\right)
\ =\ iu\sqrt{n}
-\frac{iu}{\sqrt{n}}
\Sumju \frac{1}{\bet}
-\frac{u^2}{n}
\Sumju \frac{1}{(\bet)^2}
+O\Bigg(\frac{1}{n^{3/2}}
\Sumju \frac{1}{(\bet)^3}
\Bigg)\,.
\end{equation}
We now estimate these three sums.
First, we can write
\begin{equation}
\label{eq:decompSu1}
\Bigg|\,n-\Sumju \frac{1}{\bet}\,\Bigg|
\ =\ \Bigg|\,
1-\Sumju \left(\frac{1}{\bet}-1\right)
\Bigg|
\ \leqslant\ 1
+2\Ent{\frac{n}{r_n}}
+2\sum_{j=1}^{\Ent{n/r_n}}\frac{1}{\bet}
+\sum_{j=1+\Ent{n/r_n}}^{n-\Ent{n/r_n}-1}\abs{\frac{1}{\bet}-1}\,.
\end{equation}
Following Lemma~\ref{lm:asymptotics_eigenvalues}, we have
\[\frac{1}{\bet}
\ =\ O\left(\frac{n^2}{{r_n}^2 j^2}\right)
+O(1)\,,\]
where the~$O$ is uniform over all~$n$ and~$j\leqslant\Ent{n/2}$, implying that
\[\sum_{j=1}^{\Ent{n/r_n}}\frac{1}{\bet}
\ =\ O\left(\frac{n^2}{{r_n}^2}\right)
\ =\ o\big(\sqrt{n}\big)\,.\]
Plugging this and the second estimate given by Lemma~\ref{lm:bound_eigenvalues} into~(\ref{eq:decompSu1}), we deduce that
\[
\Sumju \frac{1}{\bet}
\ =\ n+o\big(\sqrt{n}\big)\,.\]
The same decomposition of the two other sums in~(\ref{eq:DLSu}) yields
\[\Sumju \frac{1}{(\bet)^2}
\ =\ n
+O\left(\frac{n}{r_n}\right)
+O\left(\frac{n^4}{{r_n}^4}\right)
+O\left(\frac{n\ln n}{r_n}\right)
\ =\ n+o(n)\]
and
\[\Sumju \frac{1}{(\bet)^3}
\ =\ n
+O\left(\frac{n}{r_n}\right)
+O\left(\frac{n^6}{{r_n}^6}\right)
+O\left(\frac{n\ln n}{r_n}\right)
\ =\ o\big(n^{3/2}\big)\,.\]
Plugging everything into~(\ref{eq:DLSu}), we get
\begin{align*}
\varphi_n\left(\frac{u}{\sqrt{n}}\right)
&\ =\ iu\sqrt{n}
-\frac{iu\big(n+o(\sqrt{n})\big)}{\sqrt{n}}
-\frac{u^2\big(n+o(n)\big)}{n}
+o(1)
\\
&\ =\ -u^2+o(1)
\ =\ \ln\mathbb{E}\big[e^{iuX}\big]+o(1)\,,
\end{align*}
where~$X$ is a centred normal variable with variance~$2$, leading to the claimed convergence in law.
\end{proof}

\subsection{The threshold: proof of Lemma~\ref{lm:threshold}}
\label{sec:Threshold}

We now turn to the proof of Lemma~\ref{lm:intermediate}, proceeding as for the proof of Lemma~\ref{lm:long_range}.

\begin{proof}[Proof of Lemma~\ref{lm:intermediate}]
Assume that~$r_n\sim\lambda n^{3/4}$ with~$\lambda>0$, and consider the random variable~$Z_\lambda$ defined by~(\ref{eq:def_Z}).
For every~$u\in\R$, we have
\[\mathbb{E}\Big[e^{iuZ_\lambda}\Big]
\ =\ \mathbb{E}\big[e^{2iuY_0}\big]
\prod_{j=1}^\infty
\mathbb{E}\Bigg[\exp\left(-\frac{3iuY_j^2}{2\lambda^2\pi^2j^2}\right)\Bigg]^2
\ =\ \exp\left[-u^2-\sum_{j=1}^{\infty}\ln\left(1+\frac{3iu}{\pi^2\lambda^2j^2}\right)\right]\,.\]
Therefore, to obtain the claimed convergence in distribution, there only remains to prove that, for every~\smash{$u\in\R$}, 
we have
\begin{equation}
\label{eq:CVS}
\limn\,\varphi_n
\left(\frac{u}{\sqrt{n}}\right)
\ =\ -u^2-\sum_{j=1}^\infty\ln\left(1+\frac{12iu}{\pi^2\lambda^2 j^2}\right)\,,
\end{equation}
where~$\varphi_n$ is the function defined by~(\ref{eq:def_phi}).
Let~$u\in\R$.
We start by writing
\[\varphi_n
\left(\frac{u}{\sqrt{n}}\right)
\ =\ iu\sqrt{n}-\frac{n-1}{2}\ln\left(1+\frac{2iu}{\sqrt{n}}\right)
-\demi\Sumju\ln\Bigg[1+\frac{2iu}{(\sqrt{n}+2iu)}\left(\frac{1}{\bet}-1\right)\Bigg]\,.\]
Expanding the logarithm around~$1$, we get
\[iu\sqrt{n}-\frac{n-1}{2}\ln\left(1+\frac{2iu}{\sqrt{n}}\right)
\ =\ -u^2+o(1)\,.\]
Thus, to prove~(\ref{eq:CVS}), there only remains to show that
\begin{equation}
\label{eq:remains}
\limn\,
\demi\Sumju\ln\Bigg[1+\frac{2iu}{(\sqrt{n}+2iu)}\left(\frac{1}{\bet}-1\right)\Bigg]
\ =\ \sum_{j=1}^\infty\ln\left(1+\frac{3iu}{\pi^2\lambda^2 j^2}\right)\,.
\end{equation}
For every~$j\in\{1,\,\ldots,\,n-1\}$, we have~\smash{$\bet\leqslant 2$}, whence
\[\mathfrak{Re}\,\Bigg[1+\frac{2iu}{(\sqrt{n}+2iu)}\left(\frac{1}{\bet}-1\right)\Bigg]
\ =\ 1+\frac{4u^2}{n+4u^2}
\left(\frac{1}{\bet}-1\right)
\ \geqslant\ 1+\frac{4u^2}{n+4u^2}
\times\left(-\demi\right)
\ \geqslant\ \demi\,.\]
Yet, it follows from Taylor's Theorem that, for every complex numbers~$z,\,z'\in\C$ such that~$\mathfrak{Re}\,z\geqslant 1/2$ and~$\mathfrak{Re}\,z'\geqslant 1/2$, we have~\smash{$|\ln z-\ln z'|\leqslant 2\,|z-z'|$}.
Therefore, we have
\begin{align*}
\sum_{j=1+\Ent{n/r_n}}^{n-\Ent{n/r_n}-1}
\ln\Bigg[1+\frac{2iu}{(\sqrt{n}+2iu)}\left(\frac{1}{\bet}-1\right)\Bigg]
&\ =\ O\Bigg(\frac{1}{\sqrt{n}}\sum_{j=1+\Ent{n/r_n}}^{n-\Ent{n/r_n}-1}\abs{\frac{1}{\bet}-1}\Bigg)
\\
&\ =\ O\left(\frac{\ln n}{n^{1/4}}\right)\,,
\end{align*}
using the second estimate given by Lemma~\ref{lm:bound_eigenvalues}.
Similarly, we have
\begin{multline*}
\sum_{j=1}^{\Ent{n/r_n}}
\biggg\{
\ln\Bigg[1+\frac{2iu}{(\sqrt{n}+2iu)}\left(\frac{1}{\bet}-1\right)\Bigg]
-\ln\left(1+\frac{3iu}{\pi^2\lambda^2 j^2}\right)\biggg\}\\
\ =\ O\left(\frac{\Ent{n/r_n}}{\sqrt{n}}\right)
+O\Bigg(\sum_{j=1}^{\Ent{n/r_n}}
\abs{\frac{2}{(\sqrt{n}+2iu)\bet}
-\frac{3}{\pi^2\lambda^2 j^2}}\Bigg)\,.
\end{multline*}
To deal with this last term, we write
\[\sum_{j=1}^{\Ent{n/r_n}}
\abs{\frac{2}{(\sqrt{n}+2iu)\bet}
-\frac{3}{\pi^2\lambda^2 j^2}}
\ =\ \frac{1}{\lambda^2}
\sum_{j=1}^{\Ent{n/r_n}}
\abs{\frac{2r_n^2}{n^2\bet}
-\frac{3}{\pi^2 j^2}}
+o\Bigg(\sum_{j=1}^{\Ent{n/r_n}}\frac{1}{\bet\sqrt{n}}\Bigg)\,.\]
Using now the bound given by Lemma~\ref{lm:asymptotics_eigenvalues}, this becomes
\begin{align*}
\sum_{j=1}^{\Ent{n/r_n}}
\abs{\frac{2}{(\sqrt{n}+2iu)\bet}
-\frac{3}{\pi^2\lambda^2 j^2}}
&\ =\ O\Bigg(\sum_{j=1}^{\Ent{n/r_n}}\!
\frac{1}{n^{3/4} j^2}\Bigg)
+O\Bigg(\sum_{j=1}^{\Ent{n/r_n}}\!
\frac{1}{\sqrt{n}}\Bigg)
+o\Bigg(\sum_{j=1}^{\Ent{n/r_n}}\!
\frac{1}{j^2}\Bigg)\\
&\ =\ O\left(\frac{1}{n^{3/4}}\right)
+O\left(\frac{1}{n^{1/4}}\right)
+o(1)
\ =\ o(1)\,.
\end{align*}
By symmetry of the eigenvalues, we obtain~(\ref{eq:remains}), concluding the proof of~(\ref{eq:CVS}) and thus of the convergence in distribution.
\end{proof}

\subsection{Preliminaries to the saddle-point method}
\label{sec:MovingContour}

We now perform some preliminary computations which will be useful for the proofs of Lemmas~\ref{lm:finite_range} and~\ref{lm:intermediate}.

For every~$n\geqslant 4$, the characteristic function of~$n-A_n$ given by Lemma~\ref{lm:fourier} is integrable, implying by Fourier's inversion formula that the auxiliary variable~$n-A_n$ admits the following density with respect to the Lebesgue measure on~$\R$:
\[f_n\,:\,x\in\R
\ \longmapsto\ 
\frac{1}{2\pi}
\int_\R du\,\exp\big[-iux+\varphi_n(u)\big]\,.\]
Now, let~$g:\R\to\R$ be a bounded and continuous function and let~$(w_n)_{n\geqslant 4}$ be a sequence of real numbers such that~$w_n\geqslant 1$ for every~$n\geqslant 4$ (we will take~$w_n=1$ in the proof of Lemma~\ref{lm:finite_range} and~\smash{$w_n=r_n^{2/3}$} in the proof of Lemma~\ref{lm:intermediate}).
For every~$n\geqslant 4$, we have
\begin{align*}
\mathbb{E}\Bigg[g\left(\frac{n-A_n}{w_n}\right)\mathbb{1}_{\{A_n<n\}}\Bigg]
&\ =\ 
\int_0^n dx\,g\left(\frac{x}{w_n}\right)
\,f_n(x)
\\
&\ =\ 
\frac{1}{2\pi}
\int_0^n dx\,g\left(\frac{x}{w_n}\right)\,
\int_\R du\,\exp\big[-iux+\varphi_n(u)\big]
\\
&\ =\ 
\frac{w_n}{2\pi}
\int_0^{n/w_n} dx\,g(x)\,
\int_\R du\,\exp\Big[-w_niux+\varphi_n(u)\Big]
\numberthis\label{eq:beforeMovingContour}
\end{align*}
where~$\varphi_n$ is the function given by~(\ref{eq:def_phi}).
We now wish to move the integration contour of~$u$ from~$\R$ to~\smash{$-i\uc/w_n+\R$}, where~$\uc>0$ will be chosen later.
For every fixed~$x\in\R$ and~$n\geqslant 4$, the function~\smash{$f:u\mapsto\exp\big[-w_n iux+\varphi_n(u)\big]$} is analytic on the open set~$U_n$ defined by~(\ref{eq:def_un}).
\begin{figure}[ht]
\centering
\begin{tikzpicture}
\draw (-3,0) node[above] {$-M$};
\draw (3,0) node[above] {$M$};
\draw[dashed,->] (-6.5,0) -- (6.5,0);
\draw (5.5,0) node[below] {$u\in\R$};
\draw[dashed,->] (-6.5,-2) -- (6.5,-2);
\draw (5.5,-2) node[below] {$u\in -i\uc/w_n+\R$};
\draw[ultra thick,blue,->] (-3,0) -- (3,0);
\draw[ultra thick,cyan,->] (-3,0) -- (-3,-2) node[midway,left] {$\mathcal{C}_1$};
\draw[ultra thick,cyan,->] (-3,-2) -- (3,-2) node[near end,above] {$\mathcal{C}_2$};
\draw[ultra thick,cyan,->] (3,-2) -- (3,0) node[midway,right] {$\mathcal{C}_3$};
\draw (0,-2) node {$\bullet$} node[above] {$-i\uc/w_n$};
\draw (0,0) node {$\bullet$} node[above] {$0$};
\end{tikzpicture}
\caption{Cauchy's Theorem allows us to replace the integral on the segment~$[-M,M]$ by the integral along the three segments~$\mathcal{C}_1$,~$\mathcal{C}_2$ and~$\mathcal{C}_3$.}
\label{fig:contour}
\end{figure}
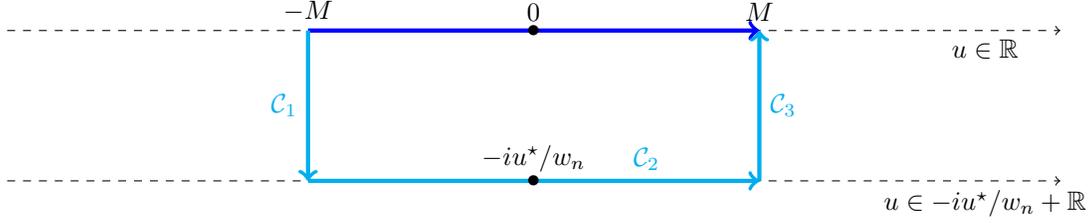
Therefore, Cauchy's theorem applied to the closed contour represented
on figure~\ref{fig:contour} ensures that, for any~$M>0$, we have
\begin{multline*}
\int_{-M}^M du\,f(u)
\ =\ 
\int_{\mathcal{C}_1}du\,f(u)
+\int_{\mathcal{C}_2}du\,f(u)
+\int_{\mathcal{C}_3}du\,f(u)
\\
\ =\ 
- \frac{i\uc}{w_n}\int_0^1 dt\,f\left(-\frac{i\uc}{w_n}\,t-M\right)
+\int_{-M}^M du\,f\left(-\frac{i\uc}{w_n}+u\right)
+ \frac{i\uc}{w_n}\int_0^1 dt\,f\left(-\frac{i\uc}{w_n}\,t+M\right)\,.
\end{multline*}
One easily checks that,~$n$ and~$x$ being fixed, we have
\[\lim_{M\to\infty}
 \int_0^1 dt\,f\left(-\frac{i\uc}{w_n}\,t-M\right)
\ =\ 
\lim_{M\to\infty}
\int_0^1 dt\,f\left(-\frac{i\uc}{w_n}\,t+M\right)
\ =\ 0\,,\]
whence
\[\int_\R du\,f(u)
\ =\ \int_\R du\,f\left(-\frac{i\uc}{w_n}+u\right)\,.\]
Going back to~(\ref{eq:beforeMovingContour}), we obtain
\begin{align*}
&\mathbb{E}\Bigg[g\left(\frac{n-A_n}{w_n}\right)\mathbb{1}_{\{A_n<n\}}\Bigg]
\\
&\qquad\ =\ 
\frac{w_n}{2\pi}
\int_0^{n/w_n} dx\,g(x)\,
\int_\R du\,\exp\Bigg[-w_niux-\uc x+\varphi_n\left(-\frac{i\uc}{w_n}+u\right)\Bigg]
\\
&\qquad\ =\ 
C_n\,
\int_0^\infty dx\,
\int_\R du\,
\mathbb{1}_{\{x<n/w_n\}}\,
g(x)
\\
&\qquad\qquad\qquad\times
\exp\Bigg[
-\frac{w_niux}{\sqrt{n}}
-\uc x
+\varphi_n\left(-\frac{i\uc}{w_n}
+\frac{u}{\sqrt{n}}\right)
-\varphi_n\left(-\frac{i\uc}{w_n}\right)
\Bigg]\,,
\numberthis\label{eq:exactExpression}
\end{align*}
where the constant~$C_n$ is given by
\[C_n
\ =\ 
\frac{w_n}{2\pi\sqrt{n}}\,
\exp\Bigg[\varphi_n\left(-\frac{i\uc}{w_n}\right)\Bigg]\,.\]
In the following two subsections, we prove the pointwise convergence of the integrand in the above formula, for a suitable choice of~$w_n$ depending on the regime of the interaction range~$r_n$.
We then want to apply the dominated convergence Theorem and, to this end, we need to check that the domination hypothesis is satisfied.
For every~$n\geqslant 4$ and every~$u\in\R$, we can write
\begin{align*}
\mathfrak{Re}\Bigg[
\varphi_n\left(-\frac{i\uc}{w_n}
+\frac{u}{\sqrt{n}}\right)
-\varphi_n\left(-\frac{i\uc}{w_n}\right)
\Bigg]
&\ =\ 
\mathfrak{Re}\Bigg[
-iu\sqrt{n}
-\demi\Sumju\ln\left(1-\frac{2iu}{\bet+2\uc/w_n}\right)\Bigg]
\\
&\ =\ 
-\frac{1}{4}\Sumju\ln\Bigg[1+\frac{4u^2}{(\bet+2\uc/w_n)^2}\Bigg]
\\
&\ \leqslant\ 
-\frac{n-1}{4}\ln\Bigg[1+\frac{4u^2}{(2+2\uc)^2}\Bigg]
\\
&\ \leqslant\ 
-\frac{3}{4}\ln\Bigg[1+\frac{u^2}{(1+\uc)^2}\Bigg]\,,
\end{align*}
where we used that~\smash{$\bet\leqslant 2$} and~$w_n\geqslant 1$ in the first inequality.
Thus, for every~$n\geqslant 4$, every~$x>0$ and every~$u\in\R$, we have
\begin{multline*}
\abs{\,
g(x)\,
\exp\Bigg[
\frac{w_niux}{\sqrt{n}}
-\uc x
+\varphi_n\left(-\frac{i\uc}{w_n}
+\frac{u}{\sqrt{n}}\right)
-\varphi_n\left(-\frac{i\uc}{w_n}\right)
\Bigg]\,
}
\\
\leqslant\ 
\sup_\R|g|\,
\exp\big(-\uc x\big)\,
\Bigg[1+\frac{u^2}{(1+\uc)^2}\Bigg]^{-3/4}\,,
\end{multline*}
which is an integrable function of~$(x,\,u)$ on~$(0,\infty)\times\R$.
Hence, there remains to study the pointwise convergence of the integrand in~(\ref{eq:exactExpression}), which is the object of the two next subsections.

\subsection{Finite range interactions: proof of Lemma~\ref{lm:finite_range}}
\label{sec:FiniteRange}

Before proving Lemma~\ref{lm:finite_range}, we state an elementary result which gives a bound on the approximation of an integral with a sum of rectangles:

\begin{lemma}[Rectangle rule]
\label{lm:rectangles}
Let~$f:[0,1]\to\C$ be a Lipschitz-continuous function with a Lipschitz constant~$K>0$. 
For every~$n\geqslant 1$, we have
\[\abs{\,\frac{1}{n}\Sumj f\left(\frac{j}{n}\right)-\int_0^1 dt\,f(t)\,}
\ \leqslant\ \frac{K}{2n}\,.\]
\end{lemma}

\begin{proof}
We write
\[\abs{\,\frac{1}{n}\Sumj f\left(\frac{j}{n}\right)-\int_0^1 dt\,f(t)\,}
\ \leqslant\ \Sumj\int_{\frac{j-1}{n}}^{\frac{j}{n}}
dt\,
\abs{\,f\left(\frac{j}{n}\right)-f(t)\,}
\ \leqslant\ \Sumj\int_{\frac{j-1}{n}}^{\frac{j}{n}} K\abs{\,\frac{j}{n}-t\,}
\ =\ \frac{K}{2n}\,,\]
which is the desired inequality.
\end{proof}

We now turn to the proof of Lemma~\ref{lm:finite_range}.

\begin{proof}[Proof of Lemma~\ref{lm:finite_range}]
Let~$r\geqslant 1$, and let~$\sigma_r>0$ be characterized by~(\ref{eq:charac_sigma}).
Note that this equation indeed has one unique solution in~$(0,\,\infty)$ because the function
\[F\,:\,x>0\ \longmapsto\ 
\int_0^1
dt\,\left(x+1-\frac{1}{r}\sum_{m=1}^r \cos(2\pi mt)\right)^{-1}
\]
is continuous and strictly decreasing and has limits
\[\lim_{x\to 0}\,F(x)
\ =\ \infty
\qquadet
\lim_{x\to \infty}\,F(x)
\ =\ 0\,,\]
so that there exists a unique value~$\sigma_r>0$ such that~\smash{$F\big(1/\sigma_r^2\big)=1$}.
The explicit value for~$\sigma_r$ in the particular case~$r=1$ follows from the formula
\[\forall a>0\quad
\forall b\in(-a,\,a)\qquad
\frac{1}{2\pi}\int_{-\pi}^{\pi}\frac{dt}{a+b\,\cos t}\ =\ \frac{1}{\sqrt{a^2-b^2}}\,,\]
which is easily obtained with a change of variable~$\theta=\tan(t/2)$, and which implies that
\[F(x)
\ =\ \frac{1}{\sqrt{(x+1)^2-1}}
\ =\ 1
\qquad\Leftrightarrow\qquad
x\ =\ \sqrt{2}-1\,.\]
Going back to the general case~$r\geqslant 1$ and taking~$w_n=1$ and~$\uc=1/(2\sigma_r^2)$, we wish to show that the function~$\varphi_n$ defined by~(\ref{eq:def_phi}) satisfies
\begin{equation}
\label{eq:CVSFiniteRange}
\forall u\in\R\qquad
\limn\,\varphi_n\left(-i\uc+\frac{u}{\sqrt{n}}\right)
-\varphi_n\big(-i\uc\big)
\ =\ -hu^2\,,
\end{equation}
for a certain constant~$h>0$.
The function~$\varphi_n$ being holomorphic on~$U_n$, it follows from Taylor's formula that, for every~$u\in\R$,
\begin{multline*}
\abs{\,
\varphi_n\left(-i\uc+\frac{u}{\sqrt{n}}\right)
-\varphi_n\big(-i\uc\big)
-\frac{u}{\sqrt{n}}\,\varphi'_n\big(-i\uc\big)
-\frac{u^2}{2n}\varphi''_n\big(-i\uc\big)
\,}
\\
\leqslant\ \frac{1}{6n^{3/2}}\sup_{t\in[0,1]}\abs{\varphi_n^{(3)}\left(-i\uc+\frac{tu}{\sqrt{n}}\right)}\,.
\numberthis\label{eq:TRI2D}
\end{multline*}
First, recalling the expression~(\ref{eq:def_phi}) of the function~$\varphi_n$ and the formula~(\ref{eq:alpha_cosine}) for~$\alp$, we write
\begin{align*}
\varphi_n'\big(-i\uc\big)
&\ =\ -in+i\Sumju\frac{1}{2\uc+\bet}
\ =\ -in-\frac{i}{2\uc}
+i\Sumj\frac{1}{2\uc+1-\alp}
\\
&\ =\ -in-\frac{i}{2\uc}
+i\Sumj f\left(\frac{j}{n}\right)\,,
\end{align*}
where the function~$f$ is given by
\[f\,:\,t\in[0,1]\ \longmapsto\ 
\left(2\uc+1-\frac{1}{r}\sum_{m=1}^r \cos(2\pi mt)\right)^{-1}\,.\]
This function being continuously differentiable on~$[0,1]$, it follows
from Lemma~\ref{lm:rectangles} that, when~$n\to\infty$,
\[\varphi_n'\big(-i\uc\big)
\ =\ -in+O(1)+in\,\int_0^1 dt\,f(t)
+O(1)
\ =\ -in+inF\left(\frac{1}{\sigma_r^2}\right)
+O(1)
\ =\ O(1)\,,\]
because~$F(1/\sigma_r^2)=1$ by definition of~$\sigma_r$.
Similarly, we have
\[\varphi_n''\big(-i\uc\big)
\ =\ -2\Sumju\frac{1}{(2\uc+\bet)^2}
\ \eqninfty\ -2hn
\]
where
\[
h\ =\ \int_0^1 dt\,
\left(2\uc+1-\frac{1}{r}\sum_{m=1}^r \cos(2\pi mt)\right)^{-2}\,,\]
following again Lemma~\ref{lm:rectangles}.
Lastly, for every~$u\in\R$, we have
\[\abs{\varphi_n^{(3)}\big(-i\uc+u\big)}
\ \leqslant\ 4\Sumju\frac{1}{(2\uc+\bet)^3}
\ \leqslant\ \frac{n}{2(\uc)^3}\,.\]
Plugging all these estimates into our Taylor expansion~(\ref{eq:TRI2D}), we obtain the claimed development~(\ref{eq:CVSFiniteRange}) of~$\varphi_n$.
Then, applying the dominated convergence Theorem to the integral
in~(\ref{eq:exactExpression}), we deduce that, for every continuous and bounded function~$g:\R\to\R$, we have
\[\limn\,\frac{1}{C_n}\mathbb{E}\Big[g(n-A_n)\mathbb{1}_{\{A_n<n\}}\Big]
\ =\ 
\int_0^\infty dx\,
\int_\R du\,
g(x)\,
e^{-\uc x
-hu^2}
\ =\ 
\sqrt{\frac{h}{\pi}}
\int_0^\infty dx\,
g(x)\,
e^{-\uc x}\,,\]
which shows the convergence in distribution of the auxiliary variable~$n-A_n$, conditioned to be positive, to an exponential variable with parameter~$\uc=1/(2\sigma_r^2)$, concluding the proof of Lemma~\ref{lm:finite_range}.
\end{proof}

\subsection{The intermediate regime: proof of Lemma~\ref{lm:intermediate}}
\label{sec:Intermediate}

Before proving Lemma~\ref{lm:intermediate}, we state an elementary result:

\begin{lemma}[An asymptotic expansion]
\label{lm:asymptotic_expansion}
When~$y\to+\infty$, we have
\[\sum_{j=1}^\infty
\frac{1}{y+j^2}
\ =\ \frac{\pi}{2\sqrt{y}}
+O\left(\frac{1}{y}\right)\,.\]
\end{lemma}

\begin{proof}
Defining, for every~$y>0$,
\[I(y)
\ =\ 
\frac{\pi}{2\sqrt{y}}
-\sum_{j=1}^\infty
\frac{1}{y+j^2}
\ =\ \int_0^\infty\frac{dt}{y+t^2}
-\sum_{j=1}^\infty
\frac{1}{y+j^2}\,,\]
we have, on the one hand,
\[I(y)
\ =\ \sum_{j=1}^\infty
\int_{j-1}^j
dt\,
\left(\frac{1}{y+t^2}-\frac{1}{y+j^2}\right)
\ \geqslant\ 0\,,\]
and, on the other hand,
\begin{align*}
I(y)
&\ \leqslant\ 
\int_0^\infty\frac{dt}{y+t^2}
-\int_0^\infty\frac{dt}{y+(t+1)^2}
\ =\ 
\int_0^\infty
\frac{(2t+1)\,dt}{(y+t^2)(y+(t+1)^2)}
\\
&\ \leqslant\ 
\int_0^\infty
\frac{(2t+1)\,dt}{(y+t^2)^2}
\ =\ \int_0^\infty
\frac{d\tau}{(y+\tau)^2}
+\int_0^\infty
\frac{dt}{(y+t^2)^2}
\\
&\ =\ 
\frac{1}{y}
\int_0^\infty
\frac{dv}{(1+v)^2}
+\frac{1}{y^{3/2}}
\int_0^\infty
\frac{dv}{(1+v^2)^2}\,,
\end{align*}
which shows that~$I(y)=O(1/y)$ when~$y\to\infty$.
\end{proof}

We now turn to the proof of Lemma~\ref{lm:intermediate}.

\begin{proof}[Proof of Lemma~\ref{lm:intermediate}]
We now assume that the interaction range is such that~$r_n/\sqrt{n}\to\infty$ and~$r_n/n^{3/4}\to 0$, we take~\smash{$w_n={r_n}^{2/3}$} and~\smash{$\uc=3^{1/3}/2^{4/3}$}, and we show that the function~$\varphi_n$ defined by~(\ref{eq:def_phi}) satisfies
\begin{equation}
\label{eq:CVSIntermediate}
\forall u\in\R\qquad
\limn\,
\varphi_n\left(
-\frac{i\uc}{{r_n}^{2/3}}
+\frac{u}{\sqrt{n}}
\right)
-\varphi_n\left(-\frac{i\uc}{{r_n}^{2/3}}\right)
\ =\ -\frac{3}{2}u^2\,.
\end{equation}
To use Taylor's formula as in~(\ref{eq:TRI2D}), we compute
\[\varphi_n'\left(-\frac{i\uc}{{r_n}^{2/3}}\right)
\ =\ -in+i
\Sumju\left(\bet+\frac{2\uc}{{r_n}^{2/3}}\right)^{-1}
\ =\ S_1+S_2\,,\]
where
\[
S_1
\ =\ 
2i
\sum_{j=1}^{\Ent{n/r_n}}
\left(\bet+\frac{2\uc}{{r_n}^{2/3}}\right)^{-1}
\qquadet
S_2
\ =\ -in+i
\sum_{j=1+\Ent{n/r_n}}^{n-1-\Ent{n/r_n}}
\left(\bet+\frac{2\uc}{{r_n}^{2/3}}\right)^{-1}
\,.\]
We now evaluate these two terms.
First, writing
\[S'_1
\ =\ 
2i
\sum_{j=1}^{\Ent{n/r_n}}
\left(\frac{2r_n^2\pi^2j^2}{3 n^2}+\frac{2\uc}{{r_n}^{2/3}}\right)^{-1}\,,\]
we have, using Lemma~\ref{lm:asymptotics_eigenvalues},
\begin{align*}
\big|S_1-S'_1\big|
&\ =\ 2
\abs{\sum_{j=1}^{\Ent{n/r_n}}
\left(\frac{2r_n^2\pi^2j^2}{3 n^2}-\bet\right)
\left(\bet+\frac{2\uc}{{r_n}^{2/3}}\right)^{-1}
\left(\frac{2r_n^2\pi^2j^2}{3 n^2}+\frac{2\uc}{{r_n}^{2/3}}\right)^{-1}}
\\
&\ \leqslant\ 
2
\sum_{j=1}^{\Ent{n/r_n}}
\abs{\frac{1}{\bet}-\frac{3 n^2}{2r_n^2\pi^2j^2}}
\ \leqslant\ 
2
\sum_{j=1}^{\Ent{n/r_n}}
\frac{K}{4}
\left(1+\frac{n^2}{{r_n}^3 j^2}\right)
\\
&\ =\ O\left(\frac{n}{r_n}\right)
+O\left(\frac{n^2}{{r_n}^3}\right)
\ =\ o\big(\sqrt{n}\big)\,.
\end{align*}
Then, defining the function
\[f\,:\,y>0\ \longmapsto\ 
\sum_{j=1}^\infty
\frac{1}{y+j^2}\]
and noting that
\[
\sum_{j=\Ent{n/r_n}+1}^\infty
\left(\frac{2r_n^2\pi^2j^2}{3 n^2}+\frac{2\uc}{{r_n}^{2/3}}\right)^{-1}
\ =\ 
O\Bigg(\frac{n^{2}}{r_n^2}
\sum_{j=\Ent{n/r_n}+1}^\infty
\frac{1}{j^2}
\Bigg)
\ =\ 
O\left(\frac{n}{r_n}\right)
\ =\ o\big(\sqrt{n}\big)\,,\]
we can write, using Lemma~\ref{lm:asymptotic_expansion},
\begin{align*}
S_1
&\ =\ S'_1+o\big(\sqrt{n}\big)
\ =\ \frac{3in^{2}}{r_n^2\pi^2}
\,f\left(\frac{3n^2\uc}{{r_n}^{8/3}\pi^2}\right)
+o\big(\sqrt{n}\big)
\\
&\ =\ \frac{3in^{2}}{r_n^2\pi^2}
\left[
\frac{\pi}{2}
\sqrt{\frac{{r_n}^{8/3}\pi^2}{3n^2\uc}}
+O\left(\frac{{r_n}^{8/3}}{n^2}\right)\right]
+o\big(\sqrt{n}\big)
\\
&\ =\ \frac{in\sqrt{3}}{2{r_n}^{2/3}\sqrt{\uc}}
+O\big({r_n}^{2/3}\big)
+o\big(\sqrt{n}\big)
\ =\ \frac{in\sqrt{3}}{2{r_n}^{2/3}\sqrt{\uc}}
+o\big(\sqrt{n}\big)\,.
\numberthis\label{eq:DLS1}
\end{align*}
We now deal with the second term.
Defining
\[
S'_2
\ =\ i
\sum_{j=1+\Ent{n/r_n}}^{n-1-\Ent{n/r_n}}\Bigg[\left(1-\alp+\frac{2\uc}{{r_n}^{2/3}}\right)^{-1}-1\Bigg]\,,\]
we have
\[S_2
\ =\ S'_2+O\left(\frac{n}{r_n}\right)
\ =\ S'_2+o\big(\sqrt{n}\big)\,.\]
Then, we can write
\begin{align*}
S'_2
&\ =\ 
i
\sum_{j=1+\Ent{n/r_n}}^{n-1-\Ent{n/r_n}}
\Bigg[
\left(1+\frac{2\uc}{{r_n}^{2/3}}\right)^{-1}
-1
+\alp\left(1+\frac{2\uc}{{r_n}^{2/3}}\right)^{-1}\left(1-\alp+\frac{2\uc}{{r_n}^{2/3}}\right)^{-1}
\Bigg]
\\
&\ =\ 
i
\sum_{j=1+\Ent{n/r_n}}^{n-1-\Ent{n/r_n}}
\Bigg[
-\frac{2\uc}{{r_n}^{2/3}}
+\frac{4(\uc)^2}{{r_n}^{4/3}}
\left(1+\frac{2\uc}{{r_n}^{2/3}}\right)^{-1}
\\
&\qquad\qquad+\alp
\left(1+\frac{2\uc}{{r_n}^{2/3}}\right)^{-1}
\left(1-\alp+\frac{2\uc}{{r_n}^{2/3}}\right)^{-1}
\Bigg]
\\
&\ =\ 
-\frac{2in\uc}{{r_n}^{2/3}}
+o\big(\sqrt{n}\big)
+O\Bigg[
\sum_{j=1+\Ent{n/r_n}}^{n-1-\Ent{n/r_n}}
\alp
\left(1-\alp+\frac{2\uc}{{r_n}^{2/3}}\right)^{-1}
\Bigg]\,.
\end{align*}
Thanks to the upper bound~(\ref{eq:upper_bound_alpha}) on~$\alp$, we know that
\[\forall j\in\Big\{1+\Ent{n/r_n},\,n-1-\Ent{n/r_n}\Big\}
\qquad
1-\alp+\frac{2\uc}{{r_n}^{2/3}}
\ \geqslant\ 
1-\alp
\ \geqslant\ 
\demi\,,\]
allowing us to expand, using again the upper bound~(\ref{eq:upper_bound_alpha}) on~$\alp$ and the first estimate of Lemma~\ref{lm:bound_eigenvalues}, and noting that~\smash{$\Sumj\alp=0$},
\begin{align*}
&\sum_{j=1+\Ent{n/r_n}}^{n-1-\Ent{n/r_n}}
\alp
\left(1-\alp+\frac{2\uc}{{r_n}^{2/3}}\right)^{-1}
\\
&\qquad\ =\ 
\sum_{j=1+\Ent{n/r_n}}^{n-1-\Ent{n/r_n}}
\alp
+O\Bigg(
\sum_{j=1+\Ent{n/r_n}}^{n-1-\Ent{n/r_n}}
\big(\alp\big)^2
\Bigg)
+O\Bigg(
\sum_{j=1}^{n-1}
\frac{\big|\alp\big|}{{r_n}^{2/3}}
\Bigg)
\\
	&\qquad\ =\ 
O\left(\frac{n}{r_n}\right)
+\sum_{j=1}^{n-1}
\alp
+O\Bigg(
\frac{n^2}{{r_n}^2}
\sum_{j=1+\Ent{n/r_n}}^{\Ent{n/2}}
\frac{1}{j^2}
\Bigg)
+O\left(\frac{n\ln n}{{r_n}^{5/3}}\right)
\ =\ 
O\left(\frac{n}{r_n}\right)\,,
\end{align*}
leading to
\[S_2
\ =\ -\frac{2in\uc}{{r_n}^{2/3}}
+o\big(\sqrt{n}\big)\,.\]
Combining this with~(\ref{eq:DLS1}), we deduce that
\[
\varphi_n'\left(-\frac{i\uc}{{r_n}^{2/3}}\right)
\ =\ \left(
\frac{\sqrt{3}}{2\sqrt{\uc}}
-2\uc\right)\frac{in}{{r_n}^{2/3}}
+o\big(\sqrt{n}\big)
\ =\ o\big(\sqrt{n}\big)\,,\]
since we took~\smash{$\uc=3^{1/3}/2^{4/3}$}.
We now study the second derivative, writing
\[\varphi_n''\left(-\frac{i\uc}{{r_n}^{2/3}}\right)
\ =\ 
-2
\Sumju\left(\bet+\frac{2\uc}{{r_n}^{2/3}}\right)^{-2}\,.
\]
Dividing the sum similarly, we show that, on the one hand,
\[-2
\sum_{j=1+\Ent{n/r_n}}^{n-1-\Ent{n/r_n}}
\left(\bet+\frac{2\uc}{{r_n}^{2/3}}\right)^{-2}
\ =\ 
-2n+O\left(\frac{n\ln n}{r_n}\right)
+O\left(\frac{n}{{r_n}^{2/3}}\right)
\ =\ -2n+o(n)\]
and, on the other hand,
\begin{align*}
-4
\sum_{j=1}^{\Ent{n/r_n}}
\left(\bet+\frac{2\uc}{{r_n}^{2/3}}\right)^{-2}
&\ =\ 
O\left(\frac{n^2}{{r_n}^2}\right)
+O\left(\frac{n}{r_n}\right)
-4
\sum_{j=1}^{\Ent{n/r_n}}
\left(\frac{2r_n^2\pi^2j^2}{3n^2}+\frac{2\uc}{{r_n}^{2/3}}\right)^{-2}
\\
&\ =\ 
o(n)
+O\big({r_n}^{4/3}\big)
+O\left(\frac{n}{r_n}\right)
-\frac{n\sqrt{3}}{\pi(\uc)^{3/2}}
\int_0^\infty
\frac{dt}{(1+t^2)^2}\\
&\ =\ 
-\frac{n\sqrt{3}}{4(\uc)^{3/2}}
+o(n)\,,
\end{align*}
whence
\[
\varphi_n''\left(-\frac{i\uc}{{r_n}^{2/3}}\right)
\ =\ 
-2n-\frac{n\sqrt{3}}{4(\uc)^{3/2}}+o(n)
\ =\ 
-3n+o(n)\,.\]
Lastly, for every~$u\in\R$, we can write
\[
\abs{\,
\varphi_n^{(3)}
\left(-\frac{i\uc}{{r_n}^{2/3}}+u\right)
\,}
\ \leqslant\ 
8\Sumju\left(\bet+\frac{2\uc}{{r_n}^{2/3}}\right)^{-3}
\ =\ S_3+S_4\,,\]
where
\[
S_3
\ =\ 
16\sum_{j=1}^{s_n}
\left(\bet+\frac{2\uc}{{r_n}^{2/3}}\right)^{-3}
\qquadet
S_4
\ =\ 
8\sum_{j=s_n+1}^{n-s_n-1}
\left(\bet+\frac{2\uc}{{r_n}^{2/3}}\right)^{-3}
\,,
\]
with~$s_n=\Ent{n/{r_n}^{4/3}}$.
One the one hand, we have
\[
S_3
\ \leqslant\ 
16\Ent{\frac{n}{{r_n}^{4/3}}}
\times\frac{{r_n}^2}{(2\uc)^3}
\ =\ O\big(n\,{r_n}^{2/3}\big)
\ =\ o\big(n^{3/2}\big)\]
and, on the other hand, using Lemma~\ref{lm:asymptotics_eigenvalues}, we can write
\begin{align*}
S_4
&\ \leqslant\ 
8\sum_{j=s_n+1}^{n-s_n-1}
\frac{1}{(\bet)^3}
\ =\ O(n)
+O\Bigg[\sum_{j=s_n+1}^\infty\left(\frac{n^2}{{r_n}^2j^2}\right)^3\Bigg]
\\
&\ =\ O(n)
+O\left(\frac{n^6}{{r_n}^6 {s_n}^5}\right)
\ =\ O\big(n\,{r_n}^{2/3}\big)
\ =\ o\big(n^{3/2}\big)\,.
\end{align*}
Thus, applying Taylor's theorem as in~(\ref{eq:TRI2D}), we obtain the pointwise convergence~(\ref{eq:CVSIntermediate}).
After this, applying the dominated convergence Theorem to the integral
in~(\ref{eq:exactExpression}), we get, for every continuous and bounded function~$g:\R\to\R$,
\begin{align*}
\limn\,
\frac{1}{C_n}\,
\mathbb{E}
\Bigg[g\left(\frac{n-A_n}{{r_n}^{2/3}}\right)\mathbb{1}_{\{A_n<n\}}\Bigg]
&\ =\ 
\int_0^\infty dx\,
\int_\R du\,
g(x)\,
e^{-\uc x
-3u^2/2}
\\
&\ =\ 
\sqrt{\frac{2\pi}{3}}
\int_0^\infty dx\,
g(x)\,
e^{-\uc x}\,,
\end{align*}
which shows the claimed convergence in law, concluding the proof of the Lemma.
\end{proof}

As explained before, the convergence in law we just proved might hold
also in the regime~$r_n=o(\sqrt{n})$.
Yet, to extend the above proof to possibly deal with this regime, a finer analysis
would be needed.
Indeed, the
hypothesis~$r_n/\sqrt{n}\to\infty$ is used many times in our
estimation of the first and
second derivative of~$\varphi_n$ (while the estimate about the third
derivative only uses~$r_n=o(n^{3/4})$).

\appendix

\section{Study of the eigenvalues}
\label{sec:study_eigenvalues}

\subsection{Proof of Lemma~\ref{lm:bound_eigenvalues}}
\label{sec:Su}

\begin{proof}[Proof of Lemma~\ref{lm:bound_eigenvalues}]
The first bound follows from the upper bound~(\ref{eq:upper_bound_alpha}) on~$\alp$, writing
\[\sum_{j=1}^n\abs{\alp}
\ \leqslant\ 
2\sum_{j=1}^{\Ent{n/2}} \abs{\alp}
\ \leqslant\ 
\frac{n}{r_n}
\sum_{j=1}^{\Ent{n/2}}\frac{1}{j}
\ =\ 
O\left(\frac{n\ln n}{r_n}\right)\,.\]
Furthermore, the upper bound~(\ref{eq:upper_bound_alpha}) on~$\alp$ ensures that when~$n/r_n<j<n-n/r_n$, we have~$\alp\leqslant 1/2$ and thus~$\bet\geqslant 1/2$, leading to
\[
\sum_{j=1+\Ent{n/r_n}}^{n-\Ent{n/r_n}-1}
\abs{\frac{1}{\bet}-1}
\ =\ 
\sum_{j=1+\Ent{n/r_n}}^{n-\Ent{n/r_n}-1}
\frac{\abs{\alp}}{\bet}
\ \leqslant\ 2
\sum_{j=1+\Ent{n/r_n}}^{n-\Ent{n/r_n}-1}
\abs{\alp}
\ \leqslant\ 
2\sum_{j=1}^n\abs{\alp}\,,\]
so that the second estimate in the statement follows from the first one.
\end{proof}

\subsection{An asymptotic formula for the eigenvalues: proof of Lemma~\ref{lm:asymptotics_eigenvalues}}

\begin{proof}[Proof of Lemma~\ref{lm:asymptotics_eigenvalues}]
Assume that~$r_n\to\infty$.
First, for every~$n$ and~$j$ such that
\begin{equation}
\label{eq:encadrJ}
\frac{n}{2r_n+1}\ <\ j\ \leqslant\ \Ent{\frac{n}{2}}\,,
\end{equation}
our formula~(\ref{eq:alpha_cos_sin_sin}) for~$\alp$ implies that
\[\abs{\alp}
\ \leqslant\ \frac{1}{r_n\,\sin\big(j\pi/n\big)}
\ \leqslant\ \frac{1}{r_n\,\sin\big(\pi/(2r_n+1)\big)}
\ \eqninfty\ \frac{2}{\pi}\ <\ 1\,,\]
using that~$r_n\to\infty$.
Therefore, uniformly for all~$n$ and~$j$ satisfying~(\ref{eq:encadrJ}), we have
\[\frac{1}{\bet}
\ =\ \frac{1}{1-\alp}
\ =\ O(1)\ =\ \frac{6n^2}{(2r_n)^2 \pi^2 j^2}+O(1)\,,\]
which implies the desired estimate.
Hence, there only remains to deal with the case~$1\leqslant j\leqslant\Ent{n/(2r_n+1)}$.
In this case, we can write, with a uniform~$O$ with respect to~$j$ and~$n$,
\[\sin\left(\frac{(2r_n+1)j\pi}{n}\right)
\ =\ \frac{(2r_n+1)j\pi}{n}
\Bigg[\,1-\frac{(2r_n+1)^2j^2\pi^2}{6n^2}+O\left(\frac{r_n^4j^4}{n^4}\right)\,\Bigg]\,.\]
Besides, since~$j\pi/n\leqslant \pi/2$, we have
\[\frac{1}{\sin(j\pi/n)}
\ =\ \frac{n}{j\pi}\Bigg[\,1+O\left(\frac{j^2}{n^2}\right)\,\Bigg]\,.\]
Using our formula~(\ref{eq:alpha_sin_sin}) for~$\alp$, we obtain that, for~$1\leqslant j\leqslant\Ent{n/(2r_n+1)}$,
\begin{align*}
\alp
\ &=\ \frac{\sin\big((2r_n+1)j\pi/n\big)}{2r_n\sin(j\pi/n)}-\frac{1}{2r_n}\\
\ &=\ \frac{(2r_n+1)j\pi/n}{2r_nj\pi/n}
\Bigg[\,1-\frac{(2r_n+1)^2j^2\pi^2}{6n^2}+O\left(\frac{r_n^4j^4}{n^4}\right)\,\Bigg]
\times\Bigg[\,1+O\left(\frac{j^2}{n^2}\right)\,\Bigg]
-\frac{1}{2r_n}\\
\ &=\ \left(1+\frac{1}{2r_n}\right)\Bigg[\,1-\frac{(2r_n+1)^2j^2\pi^2}{6n^2}+O\left(\frac{j^2}{n^2}\right)+O\left(\frac{r_n^4j^4}{n^4}\right)\,\Bigg]-\frac{1}{2r_n}\\
\ &=\ 1-\frac{(2r_n+1)^2 j^2\pi^2}{6 n^2}
+O\left(\frac{r_nj^2}{n^2}\right)
+O\left(\frac{r_n^4j^4}{n^4}\right)\\
\ &=\ 1-\frac{2r_n^2j^2\pi^2}{3n^2}
+O\left(\frac{r_nj^2}{n^2}\right)
+O\left(\frac{r_n^4j^4}{n^{4}}\right)\,,
\end{align*}
implying that
\begin{equation}
\label{eq:eqBeta}
\bet
\ =\ 1-\alp
\ =\ \frac{2r_n^2j^2\pi^2}{3n^2}
\Bigg[\,1
+O\left(\frac{1}{r_n}\right)
+O\left(\frac{r_n^2j^2}{n^2}\right)\,\Bigg]\,,
\end{equation}
with uniform~$O$ symbols for all~$1\leqslant j\leqslant\Ent{n/(2r_n+1)}$.
To obtain the result, we wish to take the inverse of this development.
To this end, we show that the quantity between brackets in~(\ref{eq:eqBeta}) is bounded away from~$0$.
Let us go back to the formula~(\ref{eq:alpha_sin_sin}) which reads
\begin{equation}
\label{eq:reecrFormuleAlp}
\alp
\ =\ \frac{\sin\big((2r_n+1)j\pi/n\big)}{2r_n\sin(j\pi/n)}-\frac{1}{2r_n}\,.
\end{equation}
A straightforward function study shows that, for every~$x\geqslant 0$,
\begin{equation}
\label{eq:encadrementSinus}
x\left(1-\frac{x^2}{6}\right)
\ \leqslant\ \sin x
\ \leqslant\ x\left(1-\frac{x^2}{6}+\frac{x^4}{120}\right)\,.
\end{equation}
From this we deduce that
\[\sin\left(\frac{j\pi}{n}\right)
\ \geqslant\ \frac{j\pi}{n}\left(1-\frac{j^2\pi^2}{6n^2}\right)\,.\]
Yet, for~$j\leqslant\Ent{n/(2r_n+1)}$, we have
\[\frac{j^2\pi^2}{6n^2}
\ \leqslant\ \frac{\pi^2}{6(2r_n+1)^2}
\ \leqslant\ \demi\]
for~$n$ large enough, because~$r_n\to\infty$.
Therefore, we get
\begin{equation}
\label{eq:majoDenom}
\sin\left(\frac{j\pi}{n}\right)^{-1}
\ \leqslant\ \frac{n}{j\pi}\left(1-\frac{j^2\pi^2}{6n^2}\right)^{-1}
\ \leqslant\ \frac{n}{j\pi}
\left(\,1+\frac{j^2\pi^2}{3n^2}\,\right)\,,
\end{equation}
using the convexity of the function~$u\mapsto(1-u)^{-1}$ over~$[0,1/2]$.
We now use the other inequality in~(\ref{eq:encadrementSinus}) to write, for~$j\leqslant\Ent{n/(2r_n+1)}$,
\begin{align*}
\sin\left(\frac{(2r_n+1)j\pi}{n}\right)
\ &\leqslant\ \frac{(2r_n+1)j\pi}{n}
\left(\,1
-\frac{(2r_n+1)^2j^2\pi^2}{6n^2}
+\frac{(2r_n+1)^4j^4\pi^4}{120n^4}\,\right)\\
\ &\leqslant\ \frac{(2r_n+1)j\pi}{n}
\Bigg[\,1-\left(\frac{\pi^2}{6}-\frac{\pi^4}{120}\right)\frac{(2r_n+1)^2j^2}{n^2}\,\Bigg]\\
\ &\leqslant\ \frac{(2r_n+1)j\pi}{n}
\left(\,1-\frac{(2r_n+1)^2j^2}{2n^2}\,\right)\,.\numberthis\label{eq:majoNum}
\end{align*}
Using~(\ref{eq:majoDenom}) and~(\ref{eq:majoNum}) in our formula~(\ref{eq:reecrFormuleAlp}) yields
\begin{align*}
\alp
\ &\leqslant\ \frac{2r_n+1}{2r_n}
\left(\,1+\frac{j^2\pi^2}{3n^2}\,\right)
\left(\,1-\frac{(2r_n+1)^2j^2}{2n^2}\,\right)-\frac{1}{2r_n}\\
\ &\leqslant\ \left(\,1+\frac{1}{2r_n}\,\right)
\left(\,1+\frac{j^2\pi^2}{3n^2}-\frac{(2r_n+1)^2j^2}{2n^2}\,\right)-\frac{1}{2r_n}\\
\ &=\ 1-\left(\,1+\frac{1}{2r_n}\,\right)
\left(\,1-\frac{2\pi^2}{3(2r_n+1)^2}\,\right)
\frac{(2r_n+1)^2j^2}{2n^2}\\
\ &=\ 1-
\left(\,1+\frac{1}{2r_n}\,\right)^3
\left(\,1-\frac{2\pi^2}{3(2r_n+1)^2}\,\right)
\frac{2r_n^2j^2}{n^2}\,.
\end{align*}
Since~$r_n\to\infty$, this implies that, for~$n$ large enough, we have
\[\forall j\leqslant \frac{n}{2r_n+1}\qquad
\bet
\ =\ 1-\alp
\ \geqslant\ \frac{r_n^2j^2}{n^2}\,.\]
This allows us to take the inverse of the development~(\ref{eq:eqBeta}), yielding, with a uniform~$O$ symbol valid for all~$1\leqslant j\leqslant\Ent{n/(2r_n+1)}$,
\[\frac{1}{\bet}
\ =\ \frac{3n^2}{2r_n^2 \pi^2 j^2}
\Bigg[\,1
+O\left(\frac{1}{r_n}\right)
+O\left(\frac{r_n^2j^2}{n^2}\right)\,\Bigg]\,,\]
which implies
\[\frac{(2r_n)^2}{n^2\bet}
\ =\ \frac{6}{\pi^2 j^2}
+O\left(\frac{1}{r_n j^2}\right)
+O\left(\frac{r_n^2}{n^2}\right)\,,\]
concluding the proof of the Lemma.
\end{proof}




\begin{thebibliography}{10}



\bibitem{ACCN88}
M.~Aizenman, J.~T. Chayes, L.~Chayes, and C.~M. Newman, \emph{Discontinuity of
  the magnetization in one-dimensional {$1/|x-y|^2$} {I}sing and {P}otts
  models}, J. Statist. Phys. \textbf{50} (1988), no.~1-2, 1--40. \MR{939480}

\bibitem{Bak96}
Per Bak, \emph{How nature works}, Copernicus, New York, 1996, The science of
  self-organized criticality. \MR{1417042}

\bibitem{BTW87}
Per Bak, Chao Tang, and Kurt Wiesenfeld, \emph{Self-organized criticality: An
  explanation of the 1/f noise}, Physical review letters \textbf{59} (1987),
  no.~4, 381.

\bibitem{Billingsley99}
Patrick Billingsley, \emph{Convergence of probability measures}, second ed.,
  Wiley Series in Probability and Statistics: Probability and Statistics, John
  Wiley \& Sons, Inc., New York, 1999, A Wiley-Interscience Publication.
  \MR{1700749}

\bibitem{BG93}
Anton Bovier and V\'{e}ronique Gayrard, \emph{The thermodynamics of the
  {C}urie-{W}eiss model with random couplings}, J. Statist. Phys. \textbf{72}
  (1993), no.~3-4, 643--664. \MR{1239568}

\bibitem{BZ97}
Anton Bovier and Milo\v{s} Zahradn\'{\i}k, \emph{The low-temperature phase of
  {K}ac-{I}sing models}, J. Statist. Phys. \textbf{87} (1997), no.~1-2,
  311--332. \MR{1453742}

\bibitem{Canning1993generalized}
A.~Canning, \emph{Generalized long-range ferromagnetic {I}sing spin models}, J.
  Phys. A \textbf{26} (1993), no.~13, 3029--3036. \MR{1236380}

\bibitem{Canning1992class}
Andrew Canning, \emph{A class of long range {I}sing spin models described by
  {C}urie-{W}eiss mean field theory}, Physica A: Statistical Mechanics and its
  Applications \textbf{185} (1992), no.~1-4, 254--260.

\bibitem{Canning92saddle}
\bysame, \emph{Saddle-point mean-field theory for long-range {I}sing spin
  models in terms of the eigenvalues and eigenvectors of the interaction
  matrix}, J. Phys. A \textbf{25} (1992), no.~18, 4723--4735. \MR{1183855}

\bibitem{CG}
Rapha\"{e}l Cerf and Matthias Gorny, \emph{A {C}urie-{W}eiss model of
  self-organized criticality}, Ann. Probab. \textbf{44} (2016), no.~1,
  444--478. \MR{3456343}

\bibitem{Copson}
E.~T. Copson, \emph{Asymptotic expansions}, Cambridge Tracts in Mathematics,
  vol.~55, Cambridge University Press, Cambridge, 2004, Reprint of the 1965
  original. \MR{2139829}

\bibitem{DM20}
Nabarun Deb and Sumit Mukherjee, \emph{Fluctuations in mean-field {I}sing
  models}, Ann. Appl. Probab. \textbf{33} (2023), no.~3, 1961--2003.
  \MR{4583662}

\bibitem{EL10}
Peter Eichelsbacher and Matthias L\"{o}we, \emph{Stein's method for dependent
  random variables occurring in statistical mechanics}, Electron. J. Probab.
  \textbf{15} (2010), no. 30, 962--988. \MR{2659754}

\bibitem{EN78}
Richard~S. Ellis and Charles~M. Newman, \emph{Limit theorems for sums of
  dependent random variables occurring in statistical mechanics}, Z. Wahrsch.
  Verw. Gebiete \textbf{44} (1978), no.~2, 117--139. \MR{0503333}

\bibitem{ENR80}
Richard~S. Ellis, Charles~M. Newman, and Jay~S. Rosen, \emph{Limit theorems for
  sums of dependent random variables occurring in statistical mechanics. {II}.
  {C}onditioning, multiple phases, and metastability}, Z. Wahrsch. Verw.
  Gebiete \textbf{51} (1980), no.~2, 153--169. \MR{566313}

\bibitem{Feller}
William Feller, \emph{An introduction to probability theory and its
  applications. {V}ol. {I}}, third ed., John Wiley \& Sons, Inc., New
  York-London-Sydney, 1968. \MR{0228020}

\bibitem{These}
Nicolas Forien, \emph{Autour de la criticit{\'e} auto-organis{\'e}e}, Ph.D.
  thesis, Universit{\'e} Paris-Saclay, 2020.

\bibitem{Frigg03}
Roman Frigg, \emph{Self-organised criticality—what it is and what it
  isn’t}, Studies in History and Philosophy of Science Part A \textbf{34}
  (2003), no.~3, 613--632.

\bibitem{GornyGaussien}
M.~Gorny, \emph{A {C}urie-{W}eiss model of self-organized criticality: the
  {G}aussian case}, Markov Process. Related Fields \textbf{20} (2014), no.~3,
  563--576. \MR{3289133}

\bibitem{GV15}
Matthias Gorny and S.~R.~S. Varadhan, \emph{Fluctuations of the self-normalized
  sum in the {C}urie-{W}eiss model of {SOC}}, J. Stat. Phys. \textbf{160}
  (2015), no.~3, 513--518. \MR{3366089}

\bibitem{KLS19a}
Zakhar Kabluchko, Matthias L\"{o}we, and Kristina Schubert, \emph{Fluctuations
  of the magnetization for {I}sing models on dense {E}rd{\H o}s-{R}\'{e}nyi
  random graphs}, J. Stat. Phys. \textbf{177} (2019), no.~1, 78--94.
  \MR{4003721}

\bibitem{KLS19b}
\bysame, \emph{Fluctuations of the magnetization for {I}sing models on {E}rd{\H
  o}s-{R}\'{e}nyi random graphs---the regimes of small {$p$} and the critical
  temperature}, J. Phys. A \textbf{53} (2020), no.~35, 355004, 37. \MR{4137544}

\bibitem{KUH63}
M.~Kac, G.~E. Uhlenbeck, and P.~C. Hemmer, \emph{On the van der {W}aals theory
  of the vapor-liquid equilibrium. {I}. {D}iscussion of a one-dimensional
  model}, J. Mathematical Phys. \textbf{4} (1963), 216--228. \MR{148416}

\bibitem{Pruessner12soc}
Gunnar Pruessner, \emph{Self-organised criticality: theory, models and
  characterisation}, Cambridge University Press, 2012.

\bibitem{SG73}
Barry Simon and Robert~B. Griffiths, \emph{The {$(\phi ^{4})_{2}$} field theory
  as a classical {I}sing model}, Comm. Math. Phys. \textbf{33} (1973),
  145--164. \MR{428998}

\bibitem{SornettePerco}
Didier Sornette, \emph{Critical phase transitions made self-organized: a
  dynamical system feedback mechanism for self-organized criticality}, Journal
  de Physique I \textbf{2} (1992), no.~11, 2065--2073.

\bibitem{WPCCJ16soc}
Nicholas~W Watkins, Gunnar Pruessner, Sandra~C Chapman, Norma~B Crosby, and
  Henrik~J Jensen, \emph{25 years of self-organized criticality: Concepts and
  controversies}, Space Science Reviews \textbf{198} (2016), no.~1-4, 3--44.


\end{thebibliography}


\providecommand{\bysame}{\leavevmode\hbox to3em{\hrulefill}\thinspace}
\providecommand{\MR}{\relax\ifhmode\unskip\space\fi MR }
\providecommand{\MRhref}[2]{%
  \href{http://www.ams.org/mathscinet-getitem?mr=#1}{#2}
}
\providecommand{\href}[2]{#2}


\begin{acks}
I wish to thank the anonymous referees for their careful reading and useful comments which helped to improve the presentation.
I also thank the German Research Foundation (project number 444084038, priority program SPP 2265) for financial support.
\end{acks}


\end{document}